\documentclass[english,11pt]{article}

\usepackage{amsfonts}
\usepackage{mathrsfs}
\usepackage{amsmath}
\usepackage{amsfonts,amssymb,color}
\usepackage{dsfont}
\usepackage{curves}
\usepackage{mathrsfs}
\usepackage{pifont}
\usepackage{graphicx}
\usepackage{float}
\usepackage{cite}
\usepackage{changes}
\usepackage{color}

\numberwithin{equation}{section}
\numberwithin{figure}{section}

\def\0{\emptyset}

\def\q{\hfill\rule{1ex}{1ex}}

\newtheorem{Theorem}{Theorem}[section]

\newtheorem{Corollary}[Theorem]{Corollary}
\newtheorem{Definition}[Theorem]{Definition}
\newtheorem{Example}[Theorem]{Example}

\newtheorem{Lemma}[Theorem]{Lemma}

\newtheorem{Proposition}[Theorem]{Proposition}

\def\emptyset{\mbox{{\rm \O}}}

\date{}

\begin{document}

\title{Unified spectral Hamiltonian results of balanced\\ bipartite  graphs and complementary graphs
\footnote{Email: liumuhuo@163.com (M. Liu), hjlai@math.wvu.edu (H.-J. Lai, Corresponding author),  wuyang850228@hotmail.com (Y. Wu)}}
\author{\small
Muhuo Liu$^{1}$,\quad Yang Wu$^{2}$, \quad   Hong-Jian Lai$^{2}$   \\\small $^{1}$ Department of   Mathematics, South China
Agricultural University,
 \\\small $^{2}$ Department of Mathematics, West
Virginia University,\\ \small Morgantown, WV, USA}
\maketitle

\begin{abstract}
There have been researches on sufficient spectral conditions for
Hamiltonian properties  and path-coverable properties  of graphs.
Utilizing the Bondy-Chv\'atal closure, we provide a unified approach
to study sufficient graph eigenvalue conditions for these properties and
sharpen former spectral results in [{\em Linear Algebra Appl.}, 432 (2010), 566-570], [{\em Linear Algebra Appl.}, 432 (2010), 2170-2173], [{\em Appl. Mech. Mater.}, 336-338 (2013), 2329-2334], [{\em Linear Algebra Appl.}, 467 (2015), 254-266],
[{\em Linear Multilinear Algebra}, 64 (2016), 2252-2269], and
[{\em J. Comb. Optim.}, 35 (2018), 1104-1127],
 among others.
\end{abstract}

\newpage

\section{Introduction}
We study simple undirected graphs, with undefined terms and notation following \cite{BoMu08}.
As in \cite{BoMu08},  $\delta(G)$, $\kappa(G)$, $\kappa'(G)$ and $\overline{G}$ denote the minimum
degree, the connectivity, the edge-connectivity and the complement of a graph $G$, respectively.
For an integer $k$, a graph $G$ is {\bf $k$-connected} (resp. {\bf $k$-edge-connected}) if
$\kappa(G) \ge k$ (resp. $\kappa'(G) \ge k$).
Throughout this paper, for an integer $s \ge 1$, let $sK_1$ be the edgeless graph with $s$ vertices.
Let $S \subseteq V(G)$ be a subset. For any vertex $u \in V(G)$, define
$N_S(u) = \{v \in S: uv \in E(G)\}$. If $H$ is a subgraph of $G$, then we use $N_H(u)$ for $N_{V(H)}(u)$.
In particular, $N_G(u) = \{v \in V(G): uv \in E(G)\}$ and $d_G(u) = |N_G(u)|$.
We often use  $N(u)$ and $d(u)$ for $N_G(u)$ and $d_G(u)$, respectively, when
$G$ is understood from the context.
A graph $G$ is {\bf nontrivial} if it has at least one edge. As in \cite{BoMu08},
$G$ is {\bf Hamiltonian} (resp., {\bf  traceable}) if $G$ contains a
spanning cycle (resp., spanning path), and is {\bf Hamilton-connected} if any
pair of distinct vertices are joined by a spanning path.

\begin{Definition} \label{def1} Let $q\geq 0$     be an integers and let $G$ be a graph.
\\
{\em(i)} $G$   is {\bf $q$-traceable} $($resp.  {\bf $q$-Hamiltonian}, {\bf $q$-Hamilton-connected}$)$ if any
removal of at most  $q$ vertices from $G$ results in a traceable graph  $($resp., a Hamiltonian
graph, a Hamilton-connected graph$)$.
\\
{\em(ii)}   $G$ is {\bf $q$-edge-Hamiltonian} if any collection of
vertex-disjoint paths with at most  $q$  edges altogether must belong to a Hamiltonian cycle in $G$.\\
{\em(iii)} $G$ is  {\bf $q$-path-coverable} if $V(G)$ can be covered by no more than $q$ vertex-disjoint paths.
\end{Definition}

Following \cite{BoMu08}, we use $G[X]$ to denote
the subgraph of $G$ induced by $X$.
By Definition \ref{def1}(i), a $q$-Hamiltonian graph is also a
$(q+1)$-traceable graph. However, a  $(q+1)$-traceable graph is not
necessarily a $q$-Hamiltonian graph. For instance, the Petersen graph
is $1$-traceable, but  not $0$-Hamiltonian. Moreover, a traceable  graph is a $0$-traceable graph,  and a Hamiltonian  graph is both  a  $0$-Hamiltonian and
a $1$-traceable  graph. If $G$ is Hamilton-connected, then for any pair of vertices $\{u,v\}$ of $G$,
there is a Hamiltonian path connecting $u$ and  $v$. Thus, $G\big[V(G)\backslash \{u,v\}\big]$
contains a  Hamiltonian path, and hence $G$ is 2-traceable.

As in \cite{BoMu08}, the {\bf join} $G\vee H$ of two disjoint graphs $G$ and $H$ is defined by
$V(G \vee H) = V(G)\cup V(H)$ and $E(G \vee H) = E(G)\cup E(H)\cup \big\{xy:$ $x\in V(G)$ and $y\in V(H)\big\}$. A {\bf $k$-regular} graph is a graph with $d_G(u)=k$ for each vertex $u\in V(G)$. For two different  nonnegative integers $p$ and $q$,
a {\bf $(p,q)$-semi-regular bipartite} graph is a bipartite graph $G$ with vertex bipartition $(U,V)$ such that
$d_G(u)=p,$ $\forall u \in U$ and $d_G(v)=q,$ $\forall v\in V$. As usual, let $K_n$, $C_n$ and $K_{k,n-k}$ be the complete graph, cycle and complete bipartite graph with  $n$ vertices, respectively.
Following \cite{MDL}, for nonnegative integers  $n$, $k$  and $s$ satisfying
$s\leq   k\leq \frac{1}{2}(n+s-2)$, define the graph $M^{k,s}_{n}$  with $n$ vertices and minimum degree $k$ as follows:
\begin{align*}
M^{k,s}_{n}= K_s \vee \big(K_{n-k-1}\cup K_{k+1-s}\big).
\end{align*}
In order to characterize the exceptional graphs in our main results, we introduce several graph families in the
following. \newpage
\begin{Definition} \label{def}
Let $n$, $k$, $p,$ $q$, $r$ be five  nonnegative integers, and $s$ be an integer.
\\
{\em(i)} Define $\mathbb{B}_{n,k,s,r}=\Big\{\overline{G_1}\vee G_2:$  $G_1=(U,V)$ is a connected
$(k-s,n-k-1)$-semi-regular bipartite graph with $n-s-1+r$ vertices and $G_2$ is a spanning subgraph of $K_{s+1-r}$,
where $0\leq r\leq s+1$ and $r\neq 1\Big\}$. In particular,
$\mathbb{B}_{n,k,-1,0}=\left\{M^{k,0}_{n}\right\}=\Big\{K_{n-k-1}\cup K_{k+1}\Big\}$.
\\
{\em(ii)}
Define $\mathbb{C}_{n,s,r}=\Big\{\overline{G_1}\vee G_2:$  $G_1$ is
a connected $(p,n-s-1-p)$-semi-regular bipartite graph  with $n-s-1+r$ vertices and $G_2$ is a
spanning subgraph of $K_{s+1-r}$, where $0\leq r\leq s+1$, $r\neq 1$  and
$1\leq p\leq \frac{n-s-1}{2}\Big\}$. In particular,
$\mathbb{C}_{n,-1,0}=\Big\{K_{p}\cup K_{n-p}:$ where $1\leq p\leq \frac{n}{2}\Big\}$.
\\
{\em(iii)} Suppose that $n=2k+1-s$ and $s\leq 1$. Define $\mathbb{H}_{n,k,s,r}=\Big\{G_1\vee G_2:$
$G_1$ is a $r$-regular graph  with $n-k+r$ vertices and $G_2$ is a spanning subgraph of
$K_{k-r}$, where $0\leq r\leq k \Big\}$.
In particular,  $\mathbb{H}_{n,k,s,k}$ is the set of all $k$-regular graphs with $n$ vertices.
\\
{\em(iv)} Let  $\mathbb{D}_{n,s,r}=\Big\{\big(\overline{G_1}\vee \overline{G_2}\big)\vee G_3:$  $G_1$ and $G_2$ are two  connected
$\frac{n-s-1}{2}$-regular graphs with $|V(G_1)|+|V(G_2)|=n-r$ and $G_3$ is a spanning subgraph of
$K_{r}$ with $\mu\big(\overline{G_3}\big)\leq n-s-1$, where $0\leq r\leq s-1\Big\}$. In particular,
$\mathbb{D}_{n,1,0}=\big\{K_{\frac{n}{2},\frac{n}{2}}\big\}$.
\\
{\em(v)}
Let  $\mathbb{W}_{n,s,r}=\Big\{\overline{G_1}\vee G_2$:   $G_1$ is a   connected  $\frac{n-s-1}{2}$-regular graph with $n-r$ vertices  and $G_2$ is a spanning subgraph of
$K_{r}$ with $\mu\big(\overline{G_2}\big)\leq n-s-1$, where $0\leq r\leq \frac{n+s-1}{2}\Big\}$. In particular,
$\mathbb{W}_{n,-1,0}$ is the set of $\left(\frac{n}{2}-1\right)$-regular graphs.
\\[2mm]
{\em(vi)} Suppose that $n > k \ge 0$, $p\geq k+1$, and let $(X,Y)$ be the vertex bipartition of
$K_{n,p+q}$ with $|X| = n$ and $|Y| = p+q$. Let $X_1 \subset X$ be a subset with $|X_1| = n-k$,
$Y_1 \subset Y$ be a subset with $|Y_1| = q\geq 1$ and $K = K_{n,p+q}[X_1 \cup Y_1]$ be the induced subgraph.
Define $B_{k,n-k;p,q} = K_{n,p+q}- E(K)$. When $k$ is understood from the context, we often write
$B_{k,n-k;p,q}$ as $Z_{p,q}$  and  define $Z^{0}_{p,q} = Z_{p,q} - e$, where $e = uv \in E(Z_{p,q})$ satisfying $d_{Z_{p,q}}(u) = n$ and $d_{Z_{p,q}}(v) = p$. To simplify the notation in the proofs, we define
\begin{equation} \label{F=Z}
F_{n,k,s} = Z_{n+s-k-1,k+1-s} \mbox{ and } F^{0}_{n,k,s}= Z^{0}_{n+s-k-1,k+1-s}.
\end{equation}
\end{Definition}
 
Following \cite{BoMu08}, we use
$G=\big[U,V\big]$ to denote a bipartite graph with vertex bipartition $(U, V)$;
and $G$ is {\bf balanced} (respectively,  {\bf almost balanced}) if  $|U|=|V|$
(respectively,  if  $|U|-|V| \in \{1, -1\}$).  
Let $p$ and $q$ be two nonnegative integers.
A bipartite graph $G=[U,V]$ is {\bf $(p,q)$-traceable} if
for any subset $S \subset G$ satisfying $|S\cap U|=p$, $|S\cap V|=q$ and $|(|U| - p) - (|V| - q)| \le 1$,
$G - S$ is traceable; and $G=[U,V]$ is {\bf $(p,q)$-Hamiltonian} if
for any subset $S \subset G$ satisfying $|S\cap U|=p$, $|S\cap V|=q$ and $|U| - p = |V| - q$,
$G - S$ is Hamiltonian.

For two graphs $G$ and $H$, we write $H \subseteq G$ if $H$ is a subgraph of $G$.
For nonnegative integers $n$ and $k$, let $\mathbb{G}_n$ be the class of graphs with $n$ vertices, and
define
the {\bf $k$-closure} of a graph $G$, denoted by $\mathscr{C}_{k}(G)$,
to be the graph obtained from $G$ by recursively
joining pairs of nonadjacent vertices whose degree sum is at least $k$ until no such pair remains nonadjacent.
By definition, $G\subseteq \mathscr{C}_{k}(G)$. A graphical property $P$ is {\bf $k$-stable}
if for any graph $G\in \mathbb{G}_n$, $G$ has Property $P$ if and only if $\mathscr{C}_{k}(G)$
has Property  $P$. It is worth noting that this definition of  $k$-stable is a slightly different from that in \cite{Bondy}.

There is also a closure concept for bipartite graphs \cite{Bondy}.
Let  $k > 0$ be an integer and $G = [U,V]$ be a bipartite graph.
The {\bf bipartite closure graph} $\mathscr{B}_k(G)$ of  $G$ is the bipartite
graph obtained from $G$ by recursively joining pairs of
nonadjacent vertices $u, v$ with $u \in U$ and $v \in V$ whose degree sum is at least $k$ until no
such pair remains nonadjacent.
By definition, $G\subseteq \mathscr{B}_{k}(G)$.

Let $A(G)$ and $D(G)$, respectively,  be the adjacency matrix and  the diagonal degree
matrix  of $G$. The   {\bf signless Laplacian matrix} of $G$ is defined to be
$Q(G)=D(G)+A(G)$. The {\bf spectral radius} of $G$, denoted by  $\rho(G)$, is the largest eigenvalue
of $A(G)$, and the {\bf signless Laplacian spectral radius} of $G$, denoted by  $\mu(G)$, is the largest eigenvalue
of $Q(G)$.
Throughout this paper, let $\alpha$ be a nonnegative real number and define $\Theta(G,\alpha)$
be the largest eigenvalue of the matrix $A(G)+\alpha D(G)$.
By definition, $\Theta(G,0)=\rho(G)$ and $\Theta(G,1)=\mu(G)$.

There have been lots of studies on graphical properties warranted by various kind of graph spectral conditions.
Our current research is motivated by these studies,
as revealed in the subsections in this section. We will
have brief literature reviews on the relationship between graphical properties and
the eigenvalues of the complement of a graph in Subsection 1.1, and those of
balanced and almost balanced bipartite graphs in Subsection 1.2. As the properties
involved are possessed by complete graphs or complete balanced bipartite graphs, and are stable under
taking the corresponding Bondy-Chv\'atal closures, we in this paper investigate the relationship
between different types of graph eigenvalues and the property when a related
Bondy-Chv\'atal closure of the graph is a complete graph or a complete balanced  bipartite
graph. Our main results, as shown in Subsections 1.1 and 1.2, present unified
conclusions that generalize several former results in a number of different problems.

\subsection{Spectral results  of complement  graphs on Hamilton problem}

There have been researches on describing the Hamiltonian properties of
a graph $G$ in terms of the eigenvalues of $\overline{G}$. The following are the related pioneer results.

\begin{Theorem} \label{w1t}
Let $G$ be a graph on $n$ vertices.
\\
 {\em(i)} {\em (Fiedler and Nikiforov  \cite{Ni3})} If $\rho\big(\overline{G}\big)\leq \sqrt{n-1}$, then $G$ is traceable unless
$G= M^{0,0}_{n}$.
\\
{\em(ii)}  {\em(Fiedler and Nikiforov   \cite{Ni3})} If $\rho\big(\overline{G}\big)\leq \sqrt{n-2}$,
then $G$ is Hamiltonian unless $G= M^{1,1}_{n}$.\\
{\em(iii)} {\em(Yu and Fan   \cite{Yu1})} If $n\geq 4$ and if $\rho\big(\overline{G}\big)\leq \sqrt{\frac{(n-2)^2}{n}}$,
then $G$ is Hamilton-connected.
\\
{\em(iv)} {\em (Li and Ning  \cite{N1})}
Suppose that  $n\geq 2k+2$ and $\delta(G)\geq k\geq 0$.
If $\rho\big(\overline{G}\big)\leq \rho\left(\overline{M^{k,0}_{n}}\right)$, then $G$
is traceable unless $G\in \mathbb{B}_{n,k,-1,0}$ or $G\in \mathbb{H}_{n,k,-1,0}$.
\\
{\em(v)} {\em (Li and Ning  \cite{N1})} Suppose that $n\geq 2k+1$ and $\delta(G)\geq k\geq 1$.
If $\rho\big(\overline{G}\big)\leq \rho\left(\overline{M^{k,1}_{n}}\right)$,
then $G$ is Hamiltonian unless $G\in \mathbb{B}_{n,k,0,0}$ or $G\in \mathbb{H}_{n,k,0,0}$.
\end{Theorem}

Extensions of some of the results stated in Theorem \ref{w1t}
have been obtained by several researchers, as seen in the theorem  below.

\begin{Theorem}\label{w3t} Let $G$ be a connected  graph on $n$ vertices.
\\
{\em(i)} {\em(Yu et al.  \cite{YFX})}
Suppose that $n\geq 2k+1$ and $\delta(G)\geq k\geq q+1 \geq 1$.
If $\rho\big(\overline{G}\big)\leq \sqrt{(k-q)(n-k-1)}$, then $G$ is $q$-Hamiltonian
and $q$-edge-Hamiltonian unless $G\in \mathbb{B}_{n,k,q,r}$ or $G\in \mathbb{H}_{n,k,q,r}$.
\\
{\em(ii)} {\em(Yu et al. \cite{YFF}, Chen and Zhang  \cite{Zhang1})}
Suppose that $n\geq 2k$ and $\delta(G)\geq k\geq 2$.
If $\rho\big(\overline{G}\big)\leq \sqrt{(k-1)(n-k-1)}$,
then $G$ is Hamilton-connected   unless $G\in \mathbb{B}_{n,k,1,0}$
or $G\in \mathbb{H}_{n,k,1,r}$.
\end{Theorem}

Analogous adjacency and signless Laplacian spectral conditions of the completeness of a graph to warrant
similar or other properties have also been investigated.
The following results come  from Theorem 3.1, Theorem 3.4 and Corollary 3.2 of Yu et al. \cite{YFX}.

\begin{Theorem} \label{w4t} {\em(Yu et al. \cite{YFX})}
Let $G$ be a graph on $n$  vertices.
\\
{\em(i)} Suppose that $n\geq 2k+1$, $\delta(G)\geq k\geq \max\{q-1,1\}$ and $q\geq 1$.
If $\rho\big(\overline{G}\big)\leq \sqrt{(k-q+2)(n-k-1)}$, then $G$ is $q$-connected and
$q$-edge-connected unless $G\in \mathbb{B}_{n,k,q-2,r}$ or $G\in \mathbb{H}_{n,k,q-2,r}$.
\\
{\em(ii)} Suppose that $n\geq 2k+q+1$, $\delta(G)\geq k\geq 1$ and $q\geq 1$.
If $\rho\big(\overline{G}\big)\leq \sqrt{(k+q)(n-k-1)}$, then $G$ is $q$-path-coverable unless
$G\in \mathbb{H}_{n,k,-q,r}$ or $G= K_{k+1}\cup K_{n-k-1}$ when $q=1$.
\end{Theorem}

\begin{Theorem} \label{w6t}
Let $G$ be a  graph with $n$ vertices.
\\
{\em(i)} {\em(Zhou \cite{Zhou1})} If $\mu\big(\overline{G}\big)\leq n$, then $G$ is traceable
unless $G\in \mathbb{C}_{n,-1,0}$ or   $G\in \mathbb{W}_{n,-1,r}$.
\\
{\em(ii)} {\em(Zhou \cite{Zhou1})}   If $\mu\big(\overline{G}\big)\leq n-1$ and $n\geq 3$,
then $G$ is Hamiltonian   unless $G\in \mathbb{C}_{n,0,0}$ or $G\in \mathbb{W}_{n,0,r}$, where $1\leq r\leq \frac{n-1}{2}$.
\\
{\em(iii)} {\em(Yu and Fan \cite{Yu1})}  If $\mu\big(\overline{G}\big)\leq n-2$ and $n\geq 6$, then
$G$ is Hamilton-connected unless $G\in \mathbb{C}_{n,1,0}$ or  $G\in \mathbb{D}_{n,1,0}$ or $G\in \mathbb{W}_{n,1,r}$, where $1\leq r\leq \frac{n}{2}$.
\end{Theorem}

It is observed that in the theorems above, all the graphical properties
warranted by the various spectral properties satisfy certain level of stability, as shown in the
result of Bondy and Chv\'atal below.

\begin{Theorem} \label{w01p}
Let $n$ and $q$ be two integers with $n\geq 3$ and $q \ge 0$. Each of the following holds for a graph on $n$ vertices.
\\
{\em(i)} {\em(Bondy and Chv\'atal \cite{Bondy})} The property that ``$G$ is $q$-connected" is $(n +q-2)$-stable.
\\
{\em(ii)} {\em(Bondy and Chv\'atal \cite{Bondy})} The property that ``$G$ is $q$-edge-connected" is $(n+q-2)$-stable.
\\
{\em(iii)} {\em(Bondy and Chv\'atal \cite{Bondy})} The property that ``$G$ is $q$-path-coverable" is $(n-q)$-stable.
\\
{\em(iv)} {\em(Bondy and Chv\'atal \cite{Bondy})} The property that ``$G$ is $q$-edge-Hamiltonian" is $(n+q)$-stable.
\\
{\em(v)} {\em(Bondy and Chv\'atal \cite{Bondy})} The property that ``$G$ is $q$-Hamiltonian connected" is $(n+q+1)$-stable.\\
{\em(vi)} {\em(Bondy and Chv\'atal \cite{Bondy})} The property that ``$G$ is $q$-Hamiltonian" is $(n+q)$-stable.
\\
{\em(vii)} {\em(Liu et al. \cite{MDL})} The property that ``$G$ is $q$-traceable" is $(n+q-1)$-stable.
\end{Theorem}

These motivate  our current study. The main result of this paper is the following.

\begin{Theorem} \label{w11t} Let $n,$ $k$ and $s$ be three integers and let $G$ be a  graph on $n$ vertices.
\\
{\em(i)} Suppose that $n\geq \max\{2k,\,2k+1-s\}$ and $\delta(G)\geq k\geq \max\big\{s,1\big\}$.
If   $\rho\big(\overline{G}\big)\leq \sqrt{(k-s)(n-k-1)}$, then either   $\mathscr{C}_{n+s}(G)= K_{n}$
or $G\in \mathbb{B}_{n,k,s,r} \cup \mathbb{H}_{n,k,s,r}$, or
both $s=k-1$ and $G=\overline{K_{1,k-1}}\vee \overline{K_{1,k-1}}$.
\\
{\em(ii)} Suppose that $n\geq 3s+2$. If   $\mu\big(\overline{G}\big)\leq n-s-1$, then either $\mathscr{C}_{n+s}(G)= K_{n}$,
or  $G\in \mathbb{C}_{n,s,r} \cup \mathbb{D}_{n,s,r} \cup \mathbb{W}_{n,s,r}$.
\end{Theorem}

Since $K_n$ is $q$-traceable for any $0\leq q\leq n$, the corollary below follows immediately from Theorem  \ref{w01p}(vii) and
Theorem \ref{w11t} with $s=q-1$.
\begin{Corollary} \label{w11c}
Let $n,$ $k$ and $q$ be three nonnegative  integers and $G$ be a graph with $|V(G)|=n$.
\\
{\em(i)} Suppose that $n \ge \max\{2k,\,2k+2-q\}$ and $\delta(G)\geq k\geq \max\big\{q-1,1\big\}$.
If   $\rho\big(\overline{G}\big)\leq \sqrt{(k+1-q)(n-k-1)}$, then $G$ is $q$-traceable,
unless $G\in \mathbb{B}_{n,k,q-1,r} \cup \mathbb{H}_{n,k,q-1,r}$ or both $G= \overline{K_{1,k-1}}\vee \overline{K_{1,k-1}}$ and $q=k$.
\\
{\em(ii)} Suppose that  $n\geq 3q-1$. If   $\mu(\overline{G})\leq n-q$, then $G$ is $q$-traceable unless
$G\in \mathbb{C}_{n,q-1,r}\cup  \mathbb{D}_{n,q-1,r}\cup \mathbb{W}_{n,q-1,r}$.
\end{Corollary}

As the complete graph has all the properties listed in Theorem \ref{w01p},
Theorem \ref{w11t} generalizes the corresponding results in Theorems \ref{w1t}, \ref{w3t},
\ref{w4t} and \ref{w6t}, when $s$ is taking different appropriate values.
Motivated by Theorem \ref{w6t},  it is natural  to consider whether
the possibility that `` $G\in \mathbb{W}_{n,q-1,r}$" can be removed from the
statement of Corollary \ref{w11c}(ii). The following
example suggests that the answer is negative.

\begin{Example}\label{21e}
Let $n$ and $q$ be two integers such that $q\geq 2$,
$n\geq 3q-5$ and $n+q$ is even.  If $G_1$ is a $(q-2)$-regular graph with  $\frac{n+q-2}{2}$ vertices,  then $G= G_1\vee \left(\frac{n-q+2}{2}\right)K_{1}$ is $\frac{n+q-2}{2}$-regular, and hence $G\in \mathbb{W}_{n,q-1,r}$. Note that any deletion of $q$ vertices from $G_1$ to $G$  results in a non-traceable graph. Thus,  $G$  is not $q$-traceable.
\end{Example}

\subsection{Spectral results of balanced bipartite graphs on Hamilton problem}

Researches on predicting traceable and Hamiltonian bipartite graphs by
graph spectral conditions have been attracted many researchers, as seen in  \cite{Lu1,LiuR,LR,N1,N3},
among others.
The following theorem displays some of the spectral results on Hamiltonian properties of balanced bipartite graphs.

\begin{Theorem} \label{01t} Let $G$ be a balanced bipartite graph  on $2n$ vertices.
\\
{\em(i)} {\em(Liu et al.  \cite{LiuR})}  If $n\geq 3$, $\delta(G)\geq 1$ and $\rho(G)\geq \sqrt{n^2-2n+3}$, then $G$ is traceable.
\\
{\em(ii)} {\em(Li and Ning  \cite{Ning2017})}  If $n\geq (k+2)^2$, $\delta(G)\geq k\geq 2$ and either  $\rho(G)\geq \rho(F_{n,k,0})$ or $\mu(G)\geq \mu(F_{n,k,0})$, then $G$ is traceable unless $G=F_{n,k,0}$.
\\
{\em(iii)} {\em(Jiang et al.  \cite{Yu2019})}  If $n\geq \max\{(k+2)^2,\frac{k^2(k+1)}{2}+k+3\}$, $\delta(G)\geq k\geq 2$ and   $\rho(G)\geq \sqrt{n(n-k-1)}$, then $G$ is traceable unless $G=F_{n,k,0}$.
\\
{\em(iv)} {\em(Liu et al.  \cite{LiuR})} If $n\geq 4$, $\delta(G)\geq 2$ and $\rho(G)\geq \sqrt{n^2-2n+4}$, then $G$ is Hamiltonian  unless $G= B_{2,n-2;n-2,2}$.
\\
{\em(v)} {\em(Li and Ning \cite{N1})} If $n\geq (k+1)^2$, $\delta(G)\geq k\geq 1$ and  $\rho(G)\geq \rho\big(Z_{n-k,k}\big)$,
then $G$ is Hamiltonian unless $G= Z_{n-k,k}$.
\\
{\em(vi)} {\em(Ge and Ning \cite{N3})} If  $n\geq k^3+2k+4$,  $\delta(G)\geq k\geq 1$ and  $\rho(G)\geq \sqrt{n(n-k)}$,
then $G$ is Hamiltonian unless $G= Z_{n-k,k}$.
\\
{\em(vii)} {\em(Jiang et al. \cite{Yu2019})}  If $n\geq \max\{(k+1)^2,\frac{k^3}{2}+k+3\}$, $\delta(G)\geq k\geq 1$ and   $\rho(G)\geq \sqrt{n(n-k)}$, then $G$ is Hamiltonian unless $G= Z_{n-k,k}$.
\\
$(viii)$ {\em(Li and Ning \cite{N1})} If $n\geq (k+1)^2$,  $\delta(G)\geq k\geq 1$ and
$\mu(G)\geq \mu\big(Z_{n-k,k}\big)$, then $G$ is Hamiltonian unless $G= Z_{n-k,k}$.
\end{Theorem}

Our current research is also motivated by the results in Theorem \ref{01t}.
The following is a useful tool.

\begin{Theorem}
\label{01p} {\em(Bondy and Chv\'atal \cite{Bondy})}
A balanced bipartite  graph $G$ with $2n$ vertices  is Hamiltonian if and only if $\mathscr{B}_{n+1}(G)$
is Hamiltonian.
\end{Theorem}

To extend results in Theorem  \ref{01t}, we need a  more general form of Theorem \ref{01p}
in our arguments, as stated in the following Proposition \ref{11p}.

\begin{Proposition}\label{11p}
Let $G$ be a balanced bipartite graph with $2n$ vertices and $q \ge 0$ be an integer,
\\
{\em(i)}  $G$ is  $(q,q)$-Hamiltonian if and only if $\mathscr{B}_{n+q+1}(G)$ is $(q,q)$-Hamiltonian.
\\
{\em(ii)} $G$ is  $(q,q)$-traceable if and only if $\mathscr{B}_{n+q+1}(G)$ is $(q,q)$-traceable.
\end{Proposition}

The main result in this subsection is to find a unified approach as a generalization of
all the former results stated in Theorem \ref{01t}, as shown in the following theorem.

\begin{Theorem}\label{12t} Let $k$ and $s$  be two  nonnegative  integers and let
$G$ be a  balanced bipartite graph with  $|V(G)|=2n\geq 8k(k+1)$ and $\delta(G)\geq k\geq \max\{s,1\}$.
If either   $\rho(G)\geq \rho\big(F^{0}_{n,k,s}\big)$ or $\mu(G)\geq \mu\big(F^{0}_{n,k,s}\big)$, then
$\mathscr{B}_{n+s}(G)= K_{n,n}$ unless $G\in \left\{F_{n,k,s},\,F^{0}_{n,k,s}\right\}$.
\end{Theorem}

Since $K_{n,n}$ is $(q,q)$-Hamiltonian for $0\leq q\leq n-2$, from  Proposition \ref{11p}(i) and Theorem \ref{12t} we deduce the following result.
\begin{Corollary}\label{01c} Let $q$ and $k$  be two  nonnegative  integers and let
$G$ be a  balanced bipartite graph with $|V(G)| = 2n\geq 8k(k+1)$.   If  $\delta(G)\geq k\geq q+1$ and if  either   $\rho(G)\geq \rho\big(F^{0}_{n,k,q+1}\big)$
or $\mu(G)\geq \mu\big(F^{0}_{n,k,q+1}\big)$,
then $G$ is  $(q,q)$-Hamiltonian  unless    $G\in \left\{F_{n,k,q+1},\,F^{0}_{n,k,q+1}\right\}$.
\end{Corollary}

\begin{Theorem}\label{02c} Let $q$ and $k$  be two  nonnegative  integers and let
$G$ be a  balanced bipartite graph with $|V(G)| = 2n\geq \max\{6k(k+1),4(k+2)^2\}$ and    $\delta(G)\geq k\geq 2$. If either     $\rho(G)\geq \rho\big(F^{0}_{n,k,0}\big)$ or   $\mu(G)\geq \mu\big(F^{0}_{n,k,0}\big)$,
then $G$ is  traceable unless   $G\in \left\{F_{n,k,0},\,F^{0}_{n,k,0}\right\}$.
\end{Theorem}
With Proposition \ref{31p} below,
Corollary \ref{01c} and Theorem \ref{02c} extend  Theorem \ref{01t} for sufficiently large $n$.

\begin{Proposition}\label{31p} Let $n$, $k$ and $s$   be three nonnegative  integers.
If  $n\geq \max\{\frac{1}{2}\big(k^2+4\big)(k+1),(k+1)(k-s+2)+2\}$ and $k\geq \max\{s,1\}$,  then
$\rho\left(F^{0}_{n,k,s}\right)<\sqrt{n(n+s-k-1)}$.
\end{Proposition}

In \cite{MDL}, studies have been done on the relationship between
Hamiltonian properties of a graph $G$ and the value of $\Theta(G,\alpha)$, the
largest eigenvalue of the matrix $A(G) + \alpha D(G)$ for a real number $\alpha$.
To further the studies in \cite{MDL}, we in this paper will show a lower  bound to
$\Theta(G,\alpha)$ that assures a balanced bipartite graph $G$ to be
$(q,q)$-traceable as well as to be $(q,q)$-Hamiltonian.
Towards this aim,   for integers $k$, $n$, $s$ and a real number $\alpha$, we define
\begin{eqnarray} \label{Theta}
\varepsilon_0(s) & = & n(n+s-k-2)+(k+1)(k+2-s),\,\,\text{and}
\\ \nonumber
\Theta_0(s) & = & \alpha \left(\frac{\varepsilon_0(s)}{n}+n\right)+(1-\alpha)\sqrt{\varepsilon_0(s)}.
\end{eqnarray}

\begin{Theorem}\label{11t} Let $k$ and $s\geq -2$  be two integers and let  $G$ be a  balanced bipartite graph
with $|V(G)| = 2n\geq 6k+8$ and  $\delta(G)\geq k\geq \max\big\{|s|,1\big\}$.
If  $\Theta(G,\alpha)>\Theta_0(s)$ and $0\leq \alpha\leq 1$,
then either $\mathscr{B}_{n+s}(G)= K_{n,n}$ or $G\subseteq F_{n,k,s}$.
\end{Theorem}

As $K_{n,n}$ is  $(q,q)$-Hamiltonian for $0\leq q\leq n-2$,
Corollary \ref{11c} follows immediately from Proposition \ref{11p}(i) and Theorem \ref{11t}.

\begin{Corollary}\label{11c} Let $q$ and $k$  be two nonnegative  integers and let
$G$ be a  balanced bipartite graph with $|V(G)| = 2n\geq 6k+8$.   If $\delta(G)\geq k\geq q+1$ and if both  $\Theta(G,\alpha)>\Theta_0(q+1) $ and $0\leq \alpha\leq 1$,
then $G$ is $(q,q)$-Hamiltonian unless  $G\subseteq F_{n,k,q+1}$.
\end{Corollary}

The following theorem summarizes some of the former results using the spectral radius $\rho(G)$ or the
signless Laplacian spectral radius $\mu(G)$ to study the traceability of an almost balanced bipartite graph $G$.

\begin{Theorem} \label{02t}
Let $G\big[U,\,V\big]$ be an almost  balanced bipartite graph with $|V(G)| = 2n-1$.
\\
{\em(i)} {\em(Liu et al.  \cite{LiuR})}
Suppose that $n\geq 4$, $\delta(G)\geq 1$, and for any $v \in V$,
$d_G(v) \ge 2$. If $\rho(G)\geq \sqrt{n^2-3n+4}$, then $G$ is
traceable unless $G\in \big\{B_{2,n-2;n-3,2},\,B_{2,n-2;n-2,1}\big\}$.
\\
{\em(ii)} {\em(Yu et al. \cite{YFF})} Suppose
that $n\geq \max\left\{\frac{1}{2}(k^3+2k+4),\,(k+1)^2\right\}$ and $\delta(G)\geq k\geq 1$. If
$\rho(G)> \sqrt{n(n-k-1)}$, then $G$ is traceable unless $G= Z_{n-k-1,k}$.
\\
{\em(iii)} {\em(Yu et al. \cite{YFF})}
Suppose that $n\geq (k+1)^2$ and $\delta(G)\geq k\geq 1$. If
$\mu(G)>2n-k-2+ \frac{(k+1)^2}{n}$, then $G$ is traceable unless $G\subseteq Z_{n-k-1,k}$.
\end{Theorem}

This also motivates our research along the same line. For a real number $\alpha$, define
\begin{eqnarray} \label{omega}
\Omega(\alpha) & =& \alpha \left(2n+q-k-2+ \frac{(k+1)(k+1-q)}{n}\right)
\\ \nonumber
& \; & +(1-\alpha)\sqrt{n(n+q-k-2)+(k+1)(k+1-q)}.
\end{eqnarray}
Our main result in this direction is Theorem \ref{13t}, which generalizes Theorem \ref{02t}
when $n$ is sufficiently large.

\begin{Theorem}\label{13t} Let  $q$ and $k$ be two nonnegative integers and let
$G$ be an almost  balanced bipartite graph with $|V(G)| = 2n-1$ and $\delta(G)\geq k\geq q+1$.
\\
{\em(i)} If $n \geq 3k+4$,  $0\leq \alpha\leq 1$, and $\Theta(G,\alpha)>\Omega(\alpha)$,
then $G$ is $(q,q)$-traceable  unless  $G\subseteq Z_{n+q-k-1,k-q}$.
\\
{\em(ii)} If $n\geq 4k(k+1)$, and if either   $\rho(G)\geq \rho\left(Z^{0}_{n+q-k-1,k-q}\right)$
or $\mu(G)\geq \mu\left(Z^{0}_{n+q-k-1,k-q}\right)$,
then $G$ is  $(q,q)$-traceable   unless  $G\in \left\{Z_{n+q-k-1,k-q},\,Z^{0}_{n+q-k-1,k-q}\right\}$.
\end{Theorem}

The organization of this paper is as follows.
In Section 2, we present the proof to Theorem  \ref{w11t}. Proposition \ref{11p} will be
justified in Section 3. Section 4 is denoted to the verification of Theorem \ref{11t}.
Utilizing Theorem \ref{11t}, we then present the proof  of Theorem \ref{12t}   in Section 5.
In Section 6, we prove Theorem \ref{02c} and  Proposition \ref{31p}, and then   we complete the proof of Theorem \ref{13t} in Section 7.

\section{The proof  of Theorem  \ref{w11t}}

We start with a few additional lemmas, which are needed in our arguments.

\begin{Lemma}  \label{w21l}
Let $G$ be a graph with $|E(G)| > 0$. Each of the following holds:
\\
{\em(i)} {\em(Cvetkovi\'c et al. \cite{DC1})}
$\mu(G) \geq \min\big\{d(u)+d(v): uv \in E(G)\big\}$.
Moreover, if $G$ is connected, then equality holds if and only if $G$ is regular or semi-regular bipartite.
\\
{\em(ii)} {\em(Li and Ning \cite{N1})}
$\rho(G) \geq \min\big\{\sqrt{d(u)d(v) }: uv \in E(G)\big\}$.
Moreover, if $G$ is connected, then equality holds if and only if $G$ is regular or semi-regular bipartite.
\end{Lemma}

\begin{Lemma}  \label{w22l} Let $p$ and $q$ be integers with $p\geq q\geq 1$.
If $G$ is a $(p,q)$-semi-regular bipartite graph, then $|V(G)| \neq q+p+1$.
\end{Lemma}

\noindent {\bf Proof.} By contradiction, let $G=[U,V]$ be a $(p,q)$-semi-regular bipartite graph
with $|V(G)| = q+p+1$. By the definition of a $(p,q)$-semi-regular bipartite graph, we have
\[
|U|+|V|=p+q+1,\,p|U|=q|V|,\,|U|\geq q\geq 1,\,\text{and}\,\,|V|\geq p\geq 1.
\]
Since  the order of a
$(p, q)$-semi-regular bipartite graph $[U, V ]$ is $p+q +1$, either $(|U|, |V |) =
(q, p + 1)$ or $(|U|, |V|) = (q + 1, p)$, but, as $p > 0$ and $q > 0$, neither of the two possibilities is
consistent with $p|U|= q|V|$.
\q

\noindent {\bf Proof of Theorem \ref{w11t}.}
For the sake of  notational simplicity, throughout the proof, we let $H = \mathscr{C}_{n+s}(G)$. Our argument is to assume that
$H \neq K_n$ to prove that in Theorem \ref{w11t} (i), $G\in \mathbb{B}_{n,k,s,r} \cup \mathbb{H}_{n,k,s,r}$ unless $G=\overline{K_{1,k-1}}\vee \overline{K_{1,k-1}}$ and $s=k-1$,
and in Theorem \ref{w11t} (ii), $G\in \mathbb{C}_{n,s,r} \cup \mathbb{D}_{n,s,r} \cup \mathbb{W}_{n,s,r}$.

Since $H\neq K_n$,  $\overline{H}$ contains at least  one non-trivial component.
We shall let $F$ denote a non-trivial component of $\overline{H}$.
For any $u, v \in V(H)$ with $uv \notin E(H)$ and
$d_H(u)\geq d_H(v)$, as $H = \mathscr{C}_{n+s}(G)$, we conclude that
$d_H(u) + d_H(v) \le n + s -1$, and so for any edge $uv\in E\big(\overline{H}\big)$,
\begin{align}\label{51e}
d_{\overline{H}}(u)+d_{\overline{H}}(v)\geq 2(n-1)-(n+s-1)=n-s-1.
\end{align}

\noindent {\bf Proof of Theorem \ref{w11t} (i).}  Our  proof of Theorem \ref{w11t} (i) takes an approach
similar to those in the justifications of Theorem 1.6(ii) in  \cite{N1} and of Theorem 3.1 in \cite{YFX}.
Here, for the completeness of the proof, we present it in detail.
By (\ref{51e}), we have
\begin{align}\label{52e}
d_{\overline{H}}(u)d_{\overline{H}}(v)\geq d_{\overline{H}}(u)\big(n-s-1-d_{\overline{H}}(u)\big).
\end{align}

Since $\delta(H) \ge \delta(G) \geq k$, we have $d_{\overline{H}}(v)\leq n-k-1$. This,
together with (\ref{51e}),
implies  that $d_{\overline{H}}(v)\geq d_{\overline{H}}(u) \geq n -s-1-d_{\overline{H}}(v)\geq  k-s$.
Hence for each $uv \in E\big(\overline{H}\big)$, we have
\[
k-s\leq d_{\overline{H}}(u)\leq d_{\overline{H}}(v)\leq  n-k-1.
\]
Let $\Phi(x) = x(n-s-1-x)$ with $k-s\leq x \leq n-k-1$. The concavity of quadratic functions implies that
\begin{equation}\label{53e}
\Phi(x)\geq \min\big\{\Phi(k-s),\,\Phi(n-k-1)\big\}=(k - s)(n-k-1).
\end{equation}
By Lemma \ref{w21l}(ii), and by (\ref{52e}) and (\ref{53e}), we have
\begin{eqnarray*}
\sqrt{(k - s)(n-k-1)} & \le & \min\left\{\sqrt{d_{\overline{H}}(u)d_{\overline{H}}(v)}:\,uv\in E(F) \right\}
\\
& \le & \rho\big(\overline{H}\big)\leq \rho\big(\overline{G}\big)\leq \sqrt{(k - s)(n-k-1)}.
\end{eqnarray*}
Claim 1 below follows from Lemma \ref{w21l}(ii).

\noindent {\bf Claim 1. }
For any nontrivial component $F$ of $\overline{H}$, each of the following holds.
\\
(i) $F$ is either regular or semi-regular bipartite.
\\
(ii) For any edge $uv\in E(F)$,
we have $d_{\overline{H}}(u)=k-s\leq n-k-1=d_{\overline{H}}(v)$.
\\
(iii) $n-k \leq |V(F)|\leq n$.

We shall complete the proof of  Theorem \ref{w11t}(i) by examining the following two cases.

\noindent
{\bf Case 1.}   $\overline{H}$  contains at least two non-trivial components.

Let $F_1$ and $F_2$ be two non-trivial components of $\overline{H}$.
By Claim 1,  each of $F_1$ and $F_2$ is either regular or  semi-regular bipartite,
and for any edge $uv \in E(\overline{H})$,
$d_{\overline{H}}(u)=k-s\leq n-k-1=d_{\overline{H}}(v)$.
By Claim 1 (iii), we have   $2(n-k)\leq n$, and so $2k\leq \max\{2k+1-s,2k\}\leq n\leq 2k$.
Thus $n = 2k$, $s\geq 1$ and $\overline{H}$ must have exactly two non-trivial components $F_1$ and $F_2$ with
$|V(F_1)| = |V(F_2)| = k$. Pick an $F_i \in \{F_1,F_2\}$.

If $F_i$ is regular, then $k-s=n-k-1=k-1$, and so $s=1$. As $|V(F_1)| = |V(F_2)| = k$ and
by Claim 1(ii), $F_1=F_2=K_k$. Since $\rho\big(\overline{G}\big)=\rho\big(\overline{H}\big)=k-1$
and $\overline{H}\subseteq\overline{G}$, we have
$\overline{G}= K_{k}\cup K_{k}$, and so $G\in \mathbb{H}_{n,k,1,0}$.

If $F_i$ is semi-regular bipartite, then $d_{\overline{H}}(v)=n-k-1=k-1=|V(F_i)|-1$,
and so $F_i=K_{1,k-1}$ and $1=d_{\overline{H}}(u)=k-s$. As $s=k-1$ and by Claim 1(ii), it follows that
$F_1=F_2=K_{1,k-1}$. Since $\rho\big(\overline{G}\big)=\rho\big(\overline{H}\big)=\sqrt{k-1}$
and $\overline{H}\subseteq\overline{G}$, we have
$\overline{G}= K_{1,k-1}\cup K_{1,k-1}$, and so $G=\overline{K_{1,k-1}}\vee \overline{K_{1,k-1}}$, as desired.

\noindent {\bf Case 2.}   $\overline{H}$  contains only one  non-trivial component.

Let $F$ denote this only nontrivial component of $\overline{H}$.
By Claim 1(i), $F$ is a regular or semi-regular bipartite graph.
Assume first that $F$ is a semi-regular bipartite graph.
By Claim 1(ii), $F$ is a connected  $(n-k-1,k-s)$-semi-regular bipartite graph, and
so for some integer $r$ with $0 \le r \le k$,
$|V(F)| = n-s-1+r$. It follows that $\overline{H} =F\cup (s+1-r)K_1$.
Since  $\rho\big(\overline{H}\big)= \rho\big(\overline{G}\big)$ and $\overline{H}$ is
a spanning subgraph of $\overline{G}$, we have
$F\cup (s+1-r)K_1 \subseteq \overline{G} \subseteq F \cup K_{s+1-r}$, and
so $\overline{F}\vee (s+1-r)K_1\subseteq G\subseteq \overline{F}\vee K_{s+1-r}$.
By Lemma \ref{w22l}, this implies that $G\in \mathbb{B}_{n,k,s,r}$.

Hence we may assume that $F$ is regular.
By Claim 1 (ii),  $k-s= n-k-1$ and so $2k+1-s=n\geq 2k$, implying $s\leq 1$.
By Claim 1 (iii), we conclude that $|V(F)| = n-k+r$, for some integer $r$ with $0\leq r \leq k$.
It follows that $\overline{H}=F\cup (k-r) K_{1}$. As
$\rho\big(\overline{H}\big)= \rho\big(\overline{G}\big)$ and $\overline{H}$ is
a spanning subgraph of $\overline{G}$, we have
$F\cup (k-r)K_1 \subseteq \overline{G} \subseteq F \cup K_{k-r}$, and so
$\overline{F}\vee (k-r)K_1 \subseteq G\subseteq \overline{F}\vee K_{k-r}$.
Since $\overline{F}$ is a $r$-regular graph with $|V(\overline{F})| = n+r-k$,
by Definition \ref{def}(ii), $G\in \mathbb{H}_{n,k,s,r}$.

This completes the proof of Theorem \ref{w11t} (i).
\q

\noindent {\bf Proof of Theorem \ref{w11t} (ii).}

By (\ref{51e}) and Lemma \ref{w21l}(i), we conclude that,
for each nontrivial component $F$ of $\overline{H}$,
$n-s-1\leq  \mu(F)\leq  \mu(\overline{H})\leq \mu(\overline{G})\leq n-s-1$,
$F$ is either a regular or a semi-regular bipartite graph,
and for any $uv\in E\big(\overline{H}\big)$,
\begin{align}\label{54e}
\mu(F)=\mu\big(\overline{H}\big)=\mu\big(\overline{G}\big)=d_{\overline{H}}(u)+
d_{\overline{H}}(v)=n-s-1.
\end{align}
Similar to the proof of Theorem \ref{w11t} (i), we justify
Theorem \ref{w11t} (ii) by a case analysis.

\noindent {\bf Claim 2. } If $\overline{H}$ has a semi-regular bipartite component, then
$\overline{H}$ has exactly one nontrivial component.

Assume that $F$ is a semi-regular bipartite component of $\overline{H}$.
By (\ref{54e}), $|V(F)| \ge n-s-1$. If $\overline{H} - V(F)$ contains a nontrivial component $F'$,
then by (\ref{54e}), $|V(F')| \ge \frac{1}{2}(n-s+1)$.
It follows from $n-s-1+\frac{n-s+1}{2}\leq n$ that $n\leq 3s+1$,
contrary to the assumption that $n \ge 3s+2$.
Hence $F$ is the unique non-trivial component of $\overline{H}$. This validates the claim.

\noindent
{\bf Case 1. } $\overline{H}$ has a semi-regular bipartite component.

We assume that $F$ is a semi-regular bipartite component of $\overline{H}$.
By Claim 2, $F$ is the only nontrivial
component of $\overline{H}$.
We may assume that $F$ is a connected  $(p,q)$-semi-regular graph with $1\leq p\leq q$, and
for some integer $r$ with $0\leq r \leq  s+1$, $|V(F)| = n-s-1+r$.
Thus $s\geq -1$, $\overline{H}=F\cup (s+1-r)K_1$, and $1\leq p\leq \frac{1}{2}(n-s-1)$ by (\ref{54e}).

Since  $\mu\big(\overline{G}\big)=\mu\big(\overline{H}\big)$ and $\overline{H}$
is a spanning subgraph of $\overline{G}$, we have
$F\cup (s+1-r)K_1 \subseteq \overline{G} \subseteq F \cup K_{s+1-r}$, and
so $\overline{F}\vee (s+1-r)K_1\subseteq G\subseteq \overline{F}\vee K_{s+1-r}$.
By Lemma \ref{w22l}, we conclude that $G\in \mathbb{C}_{n,s,r}$.
This proves Theorem \ref{w11t} (ii) if Case 1 occurs.

\noindent
{\bf Case 2.}   $\overline{H}$ does not have a semi-regular bipartite component.

By Lemma \ref{w21l}(i) and the assumption of Case 2,
every non-trivial  component of $\overline{H}$ is regular.
Let $F$ denote a component of  $\overline{H}$.
Then for any vertex $u\in V(F)$, by (\ref{54e}),
$d_{\overline{H}}(u)=\frac{1}{2}(n-s-1)$, and so $|V(F)|\geq \frac{1}{2}(n-s+1)$.

If $\overline{H}$ contains at least three non-trivial components,
then $\frac{3(n-s+1)}{2}\leq n$, implying $n\leq 3(s-1)$,
contrary to the assumption that $n \ge 3s+2$.
Hence $\overline{H}$ contains at most two nontrivial components.
Let $F'$ denote the possible
nontrivial component of $\overline{H}-V(F)$, if it exists.

We first suppose that $H$ is regular, and so $\overline{H}$ is $\frac{1}{2}(n-s-1)$-regular.
In this case, either $\overline{H}= F$ or
$\overline{H}=F\cup F'$, where $F$ and $F'$ are both connected  $\frac{1}{2}(n-s-1)$-regular.
Since $\overline{H}\subseteq \overline{G}$ and $\mu(\overline{G})=\mu(\overline{H})$,
it follows that $\overline{H}=\overline{G}$, and so $G=H\in \mathbb{W}_{n,s,0}$
for $\overline{H}= F$ or $G=H\in \mathbb{D}_{n,s,0}$ for $\overline{H}=F\cup F'$.

Hence we may assume  that $H$ is not regular, and so $\overline{H} \neq F$.
Assume first that $F$ and $F'$ are two nontrivial components of
$\overline{H}$  containing  $\frac{1}{2}(n-s+1)+r_1$ and $\frac{1}{2}(n-s+1)+r_2$ vertices, respectively.
Thus $\overline{H}= F\cup F'\cup (s-1-r_1-r_2)K_1$. Since
$\mu(\overline{G})=\mu(\overline{H})$ and since
$\overline{H}$ is a spanning subgraph of $\overline{G}$, we conclude that
$\overline{F}\vee \big(\overline{F'}\vee (s-1-r_1-r_2)K_1\big)
\subseteq G\subseteq \overline{F}\vee \big(\overline{F'}\vee  K_{s-1-r_1-r_2}\big)$,
and so $G\in \mathbb{D}_{n,s,r}$, where $1\leq r\leq s-1.$

Therefore, we may assume that    $F$ is the only non-trivial component of $\overline{H}$,
and so $\overline{H}=F\cup rK_1$, where $r = |V(G)\setminus V(F)|$.
Since $F$ is $\frac{1}{2}(n-s-1)$-regular, we have $\frac{1}{2}(n-s+1)\leq |V(F)|\leq n-1$.
Since  $\mu\big(\overline{G}\big)=\mu\big(\overline{H}\big)$ and
since $\overline{H}$ is a spanning subgraph of $\overline{G}$, it follows that
$F\cup rK_1 \subseteq \overline{G} \subseteq F \cup K_{r}$, and so
$\overline{F}\vee rK_1\subseteq G\subseteq \overline{F}\vee K_{r}$.
This implies that $G\in \mathbb{W}_{n,s,r}$, where    $1\leq r\leq \frac{1}{2}(n+s-1)$.
\q

\section{The Proof   of Proposition   \ref{11p}}

The following result initiated the study of the Bondy-Chv\'atal closure concept for balanced
bipartite graphs.

\begin{Lemma} \label{asr}
{\em(Lemma 7.3.5 of \cite{Asratian})}
Let $G=[U,V]$ be a balanced bipartite graph  with $2n$ vertices.  Let $u\in U$ and $v\in V$ be two non-adjacent  vertices with
  $d_G(u) + d_G(v)\geq n+1$. Then   $G$ is   Hamiltonian if and only if  $G +uv$ is Hamiltonian. \end{Lemma}

To prove  Proposition  \ref{11p}, it suffices to prove the following two lemmas.
\begin{Lemma}\label{11l}
Let $G=[U,V]$ be a balanced bipartite  graph with $2n$ vertices and $q$ be a
nonnegative  integer.  Let $w_1\in U$ and $w_2\in V$ be two vertices satisfying $w_1w_2 \notin E(G)$
and  $d_G(w_{1}) + d_G(w_{2})\geq n+q+1$. Then the following are equivalent.
\\
{\em(i)} $G$ is  $(q,q)$-Hamiltonian.
\\
{\em(ii)} $G'=G + w_{1}w_{2}$ is $(q,q)$-Hamiltonian.
\end{Lemma}
\noindent {\bf Proof.} As (i) implies (ii) by definition, it remains to show that (ii) implies (i).
Let $S \subset V(G)$ satisfying $|S\cap U|=|S\cap V|=q$ and $G_1=G\big[V(G)\setminus S\big]$.
We are to show that $G_1$ has a Hamilton cycle.

Since $G'$ is $(q,q)$-Hamiltonian,  $G'\big[V(G)\setminus S\big]$ contains a Hamilton cycle $C$.
If $C$ is not a Hamilton cycle of $G_1$, then $w_1w_2 \in E(C)$, and so this Hamilton cycle $C$
can be expressed as $C=w_{1}w_{2}\cdots w_{2n-2q}w_1$.
Since  $w_1\in U$ and $w_2\in V$, we observe that $|N_{S}(w_1)|+|N_S(w_2)| \leq |S| = 2q$.
As $d_G(w_{1}) + d_G(w_{2})\geq n+q+1$, we have $d_{G_1}(w_{1}) + d_{G_1}(w_{2})\geq  n-q+1$.

Note that $G_1$ is a   balanced bipartite  graph with $2(n-q)$ vertices. By Lemma \ref{asr}, $G_1$ is   Hamiltonian if and only if  $G_1 + w_1w_2$ is Hamiltonian. \q

\begin{Lemma}\label{12l} Let $G=[U,V]$ be a balanced bipartite
graph with $2n$ vertices and $q$ be a nonnegative  integer.
Let $w_1\in U$ and $w_2\in V$ be two vertices satisfying $w_1w_2 \notin E(G)$
and  $d_G(w_{1}) + d_G(w_{2})\geq n+q+1$. Then the following are equivalent.
\\
{\em(i)} $G$ is  $(q,q)$-traceable.
\\
{\em(ii)} $G'=G + w_{1}w_{2}$ is $(q,q)$-traceable.
\end{Lemma}

\noindent {\bf Proof.} By definition, we observe that (i) implies (ii), and so
it suffices to show that (ii) implies (i).
Let $S \subset V(G)$ satisfying $|S\cap U|=|S\cap V|=q$ and $G_1=G\big[V(G)\setminus S\big]$.
We are to show that $G_1$ has a Hamilton path.
Since $G'$ is $(q,q)$-traceable,  $G_1 + w_1w_2$ contains a Hamilton path $P$.

If $P$ is not a Hamilton path of $G_1$, then  $w_1w_2\in E(P)$. We suppose that
$P=u_{1}u_{2}\cdots u_{2n-2q}$, where $w_1=u_i$ and $w_2=u_{i+1}$.  As $G'[U, V]$ is bipartite with $u_i\in U$ and $u_{i+1}\in V$,
we observe that $|N_{S}(u_i)|+|N_S(u_{i+1})| \le |S| = 2q$.
By the assumption that $d_G(u_{i}) + d_G(u_{i+1})\geq n+q+1$, we conclude that
\begin{equation} \label{n-q+1}
\mbox{ $d_{G_1}(u_{i}) + d_{G_1}(u_{i+1})\geq  n-q+1$.}
\end{equation}

\par   \noindent
 {\bf Case 1. $u_1\in U$.}

Then, as $u_i \in U$, $i$ is odd. If $u_iu_{2n-2q}\in E(G)$, then $P - \{u_iu_{i+1}\} + \{u_iu_{2n-2q}\}$ is a Hamilton path
of $G_1$. Hence we may assume that $u_iu_{i+1}, u_iu_{2n-2q}\not\in E(G)$.  Similarly, we have $u_{i+1}u_1\not\in E(G)$.

\noindent  {\bf Claim 1. }
There is an index $j$ with either  $i+3 \leq j \leq 2n-2q-2$ or $2\leq j\leq i-3$, such that $u_{i}u_{j},$    $u_{i+1}u_{j+1}\in E(G_1)$.

Since $G$ is bipartite, we may suppose that  $N_{G_1}(u_{i})=\big\{u_{s_1},\,u_{s_2},\,\ldots,\,u_{s_p}\big\}$,
where  $2n-2q\not\in \big\{s_1,\,s_2,\,\ldots,\,s_p\big\}=\emptyset$ and for $t\in \big\{1,2,\ldots,p\big\}$, $s_t$ is even.
If Claim 1 fails, then $N_{G_1}(u_{i+1})\subseteq
\left\{u_{1},\,u_{3},\,u_{5},\,\ldots,\,u_{2n-2q-1}\right\}\setminus
\left\{u_{s_1+1},\,u_{s_2+1},\,\ldots,\,u_{s_p+1}\right\}$. Thus,
$d_{G_1}(u_{i+1})\leq n-q-p$ and so by (\ref{n-q+1}),
$n-q+1 \le d_{G_1}(u_i)+d_{G_1}(u_{i+1})\leq p+(n-q-p)=n-q$, a contradiction.
This  completes the proof of  Claim 1.

By Claim 1,   either for some $j$ with $2\leq j\leq i-3$, both $u_{i}u_{j}\in E(G_1)$ and $u_{i+1} u_{j+1}\in E(G_1)$, whence $u_1u_2\ldots u_{j}u_iu_{i-1}\ldots u_{j+1}u_{i+1} u_{i+2}\cdots u_{2n-2q}$ is a Hamiltonian path of $G\big[V(G)\setminus S\big]$; or for some $j$ with $i+3 \leq j \leq 2n-2q-2$, both $u_{i}u_{j}\in E(G_1)$ and $u_{i+1} u_{j+1}\in E(G_1)$, whence $u_1u_2\ldots u_iu_{j}u_{j-1}\ldots u_{i+1}  u_{j+1}u_{j+2} \ldots u_{2n-2q}$ is a Hamiltonian path of $G\big[V(G)\setminus S\big]$.
This proves that  Lemma \ref{12l} (ii) implies Lemma\ref{12l} (i) in this case.

\par   \noindent
 {\bf Case 2. $u_1\in V$.}

As $u_i \in U$, $i$ is even. We first justify the following claim.

\noindent  {\bf Claim 2. }
There is an index $j$ with either  $i+3 \leq j \leq 2n-2q-1$ or $1\leq j\leq i-3$, such that $u_{i}u_{j}$  and   $u_{i+1}u_{j+1}\in E(G_1)$.

Since $G$ is bipartite, we may suppose that  $N_{G_1}(u_{i})=\big\{u_{s_1},\,u_{s_2},\,\ldots,\,u_{s_p}\big\}$,
where $s_t$ is odd  for $t\in \big\{1,2,\ldots,p\big\}$.
If Claim 1 fails, then $N_{G_1}(u_{i+1})\subseteq
\left\{u_{2},\,u_{4},\,u_{6},\,\ldots,\,u_{2n-2q}\right\}\setminus
\left\{u_{s_1+1},\,u_{s_2+1},\,\ldots,\,u_{s_p+1}\right\}$. Thus,
by (\ref{n-q+1}),
$n-q+1 \le d_{G_1}(u_i)+d_{G_1}(u_{i+1})\leq p+(n-q-p)=n-q$, a contradiction.
This  completes the proof of  Claim 2.

By Claim 2,   either for some $j$ with $1\leq j\leq i-3$, both $u_{i}u_{j}\in E(G_1)$ and $u_{i+1} u_{j+1}\in E(G_1)$, whence $u_1u_2\ldots u_{j}u_iu_{i-1}\ldots u_{j+1}u_{i+1} u_{i+2}\cdots u_{2n-2q}$ is a Hamiltonian path of $G\big[V(G)\setminus S\big]$; or for some $j$ with $i+3 \leq j \leq 2n-2q-1$, both $u_{i}u_{j}\in E(G_1)$ and $u_{i+1} u_{j+1}\in E(G_1)$, whence $u_1u_2\ldots u_iu_{j}u_{j-1}\ldots u_{i+1}  u_{j+1}u_{j+2} \ldots u_{2n-2q}$ is a Hamiltonian path of $G\big[V(G)\setminus S\big]$.
Thus in any case, Lemma\ref{12l} holds always.  \q

\section{The Proof   of Theorem  \ref{11t}  }

Following the notation in \cite{BoMu08}, if  $A, B$ are disjoint subsets of $V(G)$, then
define $E_G[A, B] = \big\{xy \in E(G): x \in A$ and $y \in B\big\}$ and  $e_G(A,B) = \big|E_G[A, B]\big|$.
The functions $\varepsilon_0(s)$ and $\Theta_0(s)$, defined in (\ref{Theta}),
will be used in the arguments in this section.

Throughout this section, {\bf let $k$, $n$ and $s$ be integers, and unless otherwise stated, we
always assume that $G = [U,V]$ is a   balanced bipartite graph with $|V(G)| = 2n$ and $H = \mathscr{B}_{n+s}(G)$.}
By definition, we have
\begin{eqnarray} \label{H-1}
&\; & \delta(H)\geq \delta(G),\,\, |E(H)| \geq |E(G)|,
\\ \nonumber
&\mbox{ and }& \forall
u \in U,\,\, v \in V \mbox{ with } uv \notin E(H), \; d_{H}(u) + d_{H}(v)\leq  n+s-1.
\end{eqnarray}

\begin{Lemma}\label{21l}
If $n \geq 3k+4$,   $s\geq -2$,   $\delta(G)\geq k \geq \max\big\{|s|,1\big\}$,  and $|E(G)| >\varepsilon_0(s)$,
then $\mathscr{B}_{n+s}(G)= K_{n,n}$ unless $K_{n,n+s-k-1}\subseteq \mathscr{B}_{n+s}(G)$.
\end{Lemma}
\noindent{\bf Proof.}
We assume that $H \neq  K_{n,n}$ to prove that $K_{n,n+s-k-1}$ must be a subgraph of $H$.
Define
\[
U_0 = \big\{w \in U: d_G(w) \ge \frac{1}{2}(n+s)\big\}, n_U = |U_0|, V_0 = \big\{w \in V: d_G(w) \ge \frac{1}{2}(n+s)\big\}
\mbox{ and }
n_V = |V_0|.
\]
\noindent
{\bf Claim 1. }  $n_U \geq k+s+3$ and $n_V \geq k+s+3$.

By symmetry, it suffices to prove $n_U \geq k+s+3$. Direct counting yields that
\begin{equation} \label{eH-1}
|E(G)| \le |E(H)| = \sum_{v \in U} d_H(v) = \sum_{v \in U_0} d_H(v) + \sum_{v \in U-U_0}d_H(v) \le
nn_U+\frac{1}{2}(n+s-1)(n-n_U).
\end{equation}

It follows by (\ref{eH-1}) and by $|E(G)| >\varepsilon_0(s)$ that
\begin{align} \nonumber
n_U
&\geq \frac{2|E(G)|}{n+1-s}-\frac{n(n+s-1)}{n+1-s}
\\ \nonumber
&> \frac{2n(n+s-k-2)+2(k+1)(k+2-s)-n(n+s-1)}{n+1-s}
\\ \nonumber
&=\frac{n^2-(2k+3-s)n + 2(k + 1)(k - s + 2)}{n+1-s}\\&\label{eH2} = k + s + 2 + \frac{\Phi(n)}{n+1-s},
\end{align}
where $\Phi(n)=n^2-(2k+3-s)n + 2(k + 1)(k - s + 2)-(n+1-s)(k+s+2)
=n^2 -(3k+5)n + 2k^2 - ks + 5k + s^2 - s + 2$.

Since $n\geq 3k+4$, we have  $\Phi'(n)=2n-(3k+5)>0$, and so
$\Phi(n)\geq  \Phi\big(3k+4\big)=k(2k-s) + 2(k-1) + s(s-1)>0$.
It follows by (\ref{eH2}) that Claim 1 holds.

Let $p_0$ and $q_0$ be two positive integers such that $p_0\geq q_0$ and
$p_0+q_0=\max\big\{p+q,$ where $K_{p,q}\subseteq H\big\}$.
By Claim 1, we may assume that $p_0\geq q_0\geq k+s+3$.
Let $U'\subseteq U$ and $V'\subseteq V$ such that $H\big[U'\cup V'\big]= K_{p_0,q_0}$ with $|U'|=p_0$ and $|V'|=q_0$.
For any $v\in V\backslash V'$, if $v$ will be adjacent with every vertex of $U'$,
then a violation to the maximality of $p_0 + q_0$ occurs. Hence $v$ is not adjacent to at least
one vertex in $U'$. By the definition of the $(n+s)$-closure of $G$ and by symmetry, we  have
\begin{equation} \label{V'U'}
\forall v\in V\backslash V', \; d_H(v)\leq n+s-q_0-1 \mbox{ and }
\forall u\in U\backslash U', \; d_H(u)\leq n+s-p_0-1.
\end{equation}

\noindent
{\bf Claim 2.} $q_0\geq n+s-k-2$.

Assume that Claim 2 does not hold.  Then $k+s+3\leq q_0 \leq n+s-k-3$.
Define $\Phi_1(x) = x^2- (n+s -1)x + n(n+s-1)$. Since
$H$ is bipartite, and by (\ref{V'U'}), we have
\begin{align*}
|E(H)| & = \sum_{v \in V} d_H(v) = \sum_{v \in V'} d_H(v) + \sum_{v \in V - V'} d_H(v)\\& \leq nq_0+ \big(n-q_0\big)\big(n+s-q_0-1\big)
\\
&=q^2_0- (n+s -1)q_0 + n(n+s-1)=\Phi_1(q_0).
\end{align*}
  As $k+s+3\leq q_0 \leq n+s-k-3$, we have
$$
\Phi_1(q_0)\leq \max\Big\{\Phi_1\big(k+s+3\big),\,\,\Phi_1\big(n+s-k-3\big)\Big\}.
$$
Since $n\geq 3k+4$, we have both
$\varepsilon_0(s)-\Phi_1\big(n+s-k-3\big)=n-(2k- s + 4)\geq k+s\geq 0$, and
$\varepsilon_0(s)-\Phi_1\big(k+s+3\big)=(s + 2)(n-2k-5)\geq 0$.
Thus, $\varepsilon_0(s)< |E(G)| \leq |E(H)| \leq \Phi_1(q_0)\leq  \varepsilon_0(s)$, a contradiction.
This completes the proof of  Claim 2.

\noindent
{\bf Claim 3.}  $p_0+q_0\geq 2n+s-k-1$.

Assume that Claim 3 fails, and so $p_0+q_0\leq 2n+s-k-2$.
By Claim 2, we  have  $p_0\geq q_0\geq n+s-k-2$.
By (\ref{V'U'}),
\begin{align}
\nonumber
|E(H)| &\leq e_H\big(U',\,V'\big)+e_H\big(U\backslash U',\,V\big)+e_H\big(U,\,V\backslash V'\big)
\\&\label{21e}
\leq p_0q_0+(n+s-1-p_0)(n-p_0)+(n+s-1-q_0)(n-q_0).
\end{align}

If  $p_0\geq n+s-k$, then as $\delta(G)\geq k$, it
follows from the definition of the $(n+s)$-closure of $G$ that
each vertex of $V'$ must be adjacent to every vertex of $U$, and so $p_0=n$
and $q_0=n+s-k-2$. It follows  by  (\ref{21e}) that
$\varepsilon_0(s) < |E(G)| \leq |E(H)| \leq  n(n+s-k-2)+(k+1)(k+2-s)=\varepsilon_0(s)$, a contradiction.
Hence we may assume that $n+s-k-2\leq q_0\leq p_0\leq n+s-k-1$.

If  $p_0=q_0=n+s-k-1$, then  by (\ref{21e}) we have $|E(H)| \leq (n+s-k-1)^2+2k(k+1-s)$.
As $n\geq 3k+4$,  this leads to
\begin{align*}
\varepsilon_0(s) - |E(H)| & \ge \varepsilon_0(s)-\big((n+s-k-1)^2+2k(k+1-s)\big)
\\
& =(k-s)(n+s-2k-1)+1 \ge (k-s)((3k+4) +s-2k-1)+1
\\
& =(k - s)(k + s + 3) + 1>0.
\end{align*}
Hence $|E(G)| \leq |E(H)| <\varepsilon_0(s)$, contrary to the assumption of the lemma.

If  $p_0=n+s-k-1$ and  $q_0=n+s-k-2$, then  by (\ref{21e}) and $n\geq 3k+4$  we have
$|E(G)| \leq (n+s-k-1)(n+s-k-2)+k(k+1-s)+(k+1)(k+2-s)<\varepsilon_0(s)$,
again a contradiction.

If  $p_0=n+s-k-2=q_0$,  then  by (\ref{21e}) and $n\geq 3k+4$ we have
$|E(G)| \leq (n+s-k-2)^2+2(k+1)(k+2-s)<\varepsilon_0(s)$,
contrary to the assumption of the lemma, and so Claim 3 is justified.

If $p_0=n$, then the lemma follows from Claim 3. Assume
that $p_0\leq n-1$, and so $q_0\geq n+s-k$ by Claim 3. As $\delta(G)\geq k$,
we conclude that every vertex of $U'$ must be adjacent to
all vertices of $V$, implying that $p_0 \geq q_0=n$, contrary to
the assumption that $H \neq K_{n,n}$. \q

\begin{Theorem} \label{21t}
If  $n \geq 3k+4$, $s\geq -2$, $\delta(G)\geq k\geq \max\big\{|s|,1\big\}$ and $|E(G)| >\varepsilon_0(s)$,
then $\mathscr{B}_{n+s}(G)$ is isomorphic to a member in $\{ K_{n,n}, F_{n,k,s}\}$.
\end{Theorem}
\noindent{\bf Proof.} We assume that $H\neq K_{n,n}$ to
show that $H= F_{n,k,s}$.
Let $t$ be the largest integer such that $K_{n,t} \subseteq H$.
By Lemma  \ref{21l}, $n +s- k-1\leq t< n$.
Let $V'\subset V$ be the vertex  sets of $H$  such that $H\big[U\cup V'\big]= K_{n,t}$.
If $t\geq n+s-k$, since  every vertex in $U$ has
degree at least $t\geq n+s-k$ in $H$ and $\delta(H)\geq k$,
we have $H= K_{n,n}$, contrary to the assumption. Hence we must have $t=n +s-k-1$.

Define $U_0 =\big\{u \in U: d_H(u) \ge n +s-k\big\}$.
Since $\delta(H)\geq k$ and since every vertex in $U$ has degree at least  $n +s-k-1$ in $H$,
it follows from the definition of the $(n+s)$-closure of $G$ that
every vertex in $V\setminus V'$ has degree
exactly $k$ in $H$, and is adjacent to every   vertex in $U_0$.
This implies  that $|U_0| = k$, and so $H= F_{n,k,s}$.
\q

We need the following two lemmas to complete the proof of  Theorem \ref{11t}.

\begin{Lemma}\label{25l} {\em(Li and Ning \cite{N1})} If  $G$ is a balanced bipartite
graph with $|V(G)| = 2n$, then $\mu(G)\leq \frac{|E(G)|}{n}+n.$
\end{Lemma}

When $|V(G)| \ge 2$, let $\rho_1(G)$ and $\rho_2(G)$ denote the largest and the second largest
eigenvalues of $A(G)$, respectively.  Thus,  $\rho_1(G)=\rho(G)$.

\begin{Lemma}\label{26l} {\em(Lai, Liu and Zhou \cite{Lai2})}
If  $G$ is a   bipartite graph with $|V(G)| \ge 2$, then
$\big(\rho_1(G)\big)^{2}+\big(\rho_2(G)\big)^{2}\leq |E(G)|$.
\end{Lemma}

When $0\leq \alpha\leq 1$,  since $A(G)+\alpha D(G)=\alpha Q(G)+(1-\alpha)A(G)$,   from
the properties of Rayleigh quotients we have $\Theta(G,\alpha) \leq \alpha \mu(G)+(1-\alpha) \rho(G)$.
Thus, the corollary below follows immediately from
  Lemmas \ref{25l} and
\ref{26l}.

\begin{Corollary} \label{C4}
Let $\alpha$ be a real number with $0\leq \alpha \leq 1$.
If $|V(G)|=2n\geq 2$, then
\begin{align}\label{22e}
\Theta(G,\alpha) \leq \alpha \mu(G)+(1-\alpha) \rho(G)\leq
\alpha \left(\frac{|E(G)|}{n}+n\right)+(1-\alpha) \sqrt{|E(G)|}.
\end{align}
\end{Corollary}

Recall that $\Theta_0(s)$ and $\varepsilon_0(s)$ have been defined in (\ref{Theta}).
If $|E(G)|\leq \varepsilon_0(s)$, then Corollary \ref{C4} implies that  $\Theta(G,\alpha)\leq \Theta_0(s)$. This deduces    the following result.
\begin{Corollary} \label{22t} Let $\alpha$ be a real number with $0\leq \alpha \leq 1$.
If $|V(G)|=2n\geq 2$  and  $\Theta(G,\alpha)>\Theta_0(s)$,
then $|E(G)| > \varepsilon_0(s)$.
\end{Corollary}

\par\bigskip \noindent
{\bf Proof of Theorem \ref{11t}.}  In our hypotheses, $\Theta(G,\alpha) > \Theta_0(s)$,
hence, Theorem    \ref{21t}  and Corollary  \ref{22t}  imply that
$\mathscr{B}_{n+s}(G)$ is isomorphic to a member in $\big\{F_{n,k,s},\, K_{n,n}\big\}$.
\q

\section{The Proof  of Theorem    \ref{12t} }

Given two distinct vertices $u, v$ in a graph $G$, if $N_G(v)\setminus \big(N_{G}(u)\cup \{u\}\big)\neq \emptyset\neq N_G(u)\setminus \big(N_{G}(v)\cup \{v\}\big)$, then we construct  a new graph $G'= G'(u, v)$ by replacing
all edges $vw$ by $uw$ for each $w\in N_{G}(v)\setminus \big(N_{G}(u)\cup \{u\}\big)$. This operation is called the {\bf Kelmans
transformation} from $v$ to $u$ (See \cite{P1}).
\begin{Lemma}\label{31l} {\em(Liu et al. \cite{MDL})}
Let $G$ be a connected  graph. If  $G'$
is  a graph obtained from $G$ by some Kelmans transformation and $\alpha\geq 0$, then $\Theta(G',\alpha)>\Theta(G,\alpha).$
\end{Lemma}

In the discussion of Lemma \ref{32l} below, the notation in Definition \ref{def} (vi) will be adopted.

\begin{Lemma}\label{32l}
Let $G$   be a  graph obtained
from $Z_{p,q}$ by deleting one edge. If $p\geq
k+1$, $q\geq 1$, $\alpha \geq 0$ and $\delta(G)\geq k\geq 1$, then  $
\Theta(G,\alpha)\leq \Theta\left(Z^{0}_{p,q},\alpha\right)$, with equality
if and only if $G= Z^{0}_{p,q}$.
\end{Lemma}
\noindent{\bf  Proof.} Let $G' = Z_{p,q}$ and $G_0 = Z^{0}_{p,q}$.  Let
$e=w_{0}z_{0} \in E(G')$,  and $G = G' - e$.
It suffices to show that if $G \neq G_0$, then
\begin{align}\label{31e}
\Theta(G,\alpha)<\Theta(G_0,\alpha).
\end{align}
Let $U$ and $V$ be the bipartition of $G'$ such that $V$ contains $q$ vertices of
degree $k$ and $U$ contains  $k$ vertices of degree $p+q$ in $G'$.
Let $U'$ and  $V'$  be the vertices of degrees $p+q$ and $n$, respectively,    in $U$ and $V$   of $G'$.
Since every vertex of $V\backslash V'$ has degree $k$ and since $G\neq G_0$,
by symmetry, we may assume that $w_0\in U'$ and $z_0\in V'$.

Choose $v\in U\setminus U'$. Then, $N_G(v)\setminus \big(N_G(w_0)\cup \{w_0\}\big)=\{z_0\}$ and $N_G(w_0)\setminus \big(N_G(v)\cup \{v\}\big)\neq \emptyset$.
It is routine to verify that $G_0$ is isomorphic to the graph obtained from
$G$ by a Kelmans transformation from $v$ to $w_0$. By Lemma \ref{31l},
$\Theta(G,\alpha)< \Theta(G_0,\alpha)$, and so (\ref{31e}) holds.
\q

Let $G$ be a connected graph. For any  real number $\alpha \geq 0$,   it is well
known that $A(G)$ is nonnegative and irreducible if and only if $G$ is
connected, and thus  $A(G)+\alpha D(G)$ is a nonnegative
irreducible matrix. This implies the existence of   a unique positive unit eigenvector
$f=\big(f(v_1),f(v_2),\ldots,f(v_n)\big)^{T}$ corresponding to
$\Theta(G,\alpha)$. This vector $f$ is often called the {\bf Perron vector} of $G$.

\begin{Lemma}\label{33l}
For any integers $n, q$ and a real number $\alpha$, define a polynomial in $\theta$ as follows:
\begin{eqnarray*}
\Psi(\theta) & = & \theta^4 -2(n + q - 1)\alpha \theta^3 +\Big(\alpha^2\big(n^2 + 4nq - 3n + q^2 - 3q + 1\big)-nq+1\Big)\theta^2
\\
& \; &-\alpha\Big(\alpha^2(2nq - q - n)(n + q - 1)-nq(n +q-2)\Big)\theta +\big(\alpha^{2}-1\big)(n-1)(q-1)\big(nq\alpha^2-1\big).
\end{eqnarray*}
If   $2\leq q \leq n$ and $\alpha \geq 0$, then $\Theta(K_{n,q}-e, \alpha)$ is   the maximum root of
$\Psi(\theta)$.
\end{Lemma}
\noindent{\bf  Proof:} Denote $G=K_{n,q}-e$ and $\Theta=\Theta(K_{n,q}-e,\alpha)$.
Let $f$ be the Perron vector of $G$, and let   $U$ and $V$ be the two partite sets
of $G$ such that $|U|=n$ and $|V|=q$.   For convenience, we suppose that $e=w_0z_0$ with $w_0\in U$ and $z_0\in V$.

Let $x_1=f(w)$ for $w\in U\setminus \{w_0\}$,
let $x_2=f(w)$ for $w\in V\setminus \{z_0\}$, let   $x_3=f(w_0)$ and $x_4=f(z_0)$.
it follows from $\big(A(G)+\alpha D(G)\big)f=\Theta f$ that
\begin{align}\label{031e}\left\{
\begin{array}{llll}
 \big(\Theta-q\alpha\big)\,x_1=(q-1)x_2+x_4, \\[1mm]
\big(\Theta-n\alpha\big) \,x_2=(n-1)x_1+x_3,  \\[1mm]
\big(\Theta-(q-1)\alpha\big)\,x_3=(q-1)x_2, \\[1mm]
\big(\Theta-(n-1)\alpha\big)\,x_4=(n-1)x_1.\end{array}
\right.\end{align}

 By multiplying $\Theta-(n-1)\alpha$ in both side of the first   equation of (\ref{031e}), and then multiplying  $\Theta-(q-1)\alpha$ in both side of the second  equation of (\ref{031e}), it follows that \begin{align}\label{0031e}\left\{
\begin{array}{llll}
 \big(\Theta-q\alpha\big)\big(\Theta-(n-1)\alpha\big)\,x_1
 =(q-1)\big(\Theta-(n-1)\alpha\big)x_2+\big(\Theta-(n-1)\alpha\big)x_4, \\[1mm]
\big(\Theta-n\alpha\big)\big(\Theta-(q-1)\alpha\big) \,x_2=\big(\Theta-(q-1)\alpha\big)(n-1)x_1+\big(\Theta-(q-1)\alpha\big)x_3.\end{array}
\right.\end{align}

  By substituting the last two equation of (\ref{031e}) into   (\ref{0031e}),  we have
 \begin{align}\label{032e}\left\{\begin{array}{llll}
 \big((\Theta-q\alpha)(\Theta-(n-1)\alpha)-(n-1)\big)\,x_1
 =(q-1)(\Theta-(n-1)\alpha)x_2, \\[1mm]
\big((\Theta-n\alpha)(\Theta-(q-1)\alpha)-(q-1)\big) \,x_2=(\Theta-(q-1)\alpha)(n-1)x_1.\end{array}
\right.\end{align}

Now, by (\ref{032e}), $\Theta$ is equal to the maximum root of
$\Psi(\theta)$, as required.  \q

\begin{Corollary}\label{31c}
Let $k$ and $s$ be two nonnegative integers such that $k\geq \max\{s,1\}$.
Each of the following holds.
\\
{\em(i)} If   $n\geq (k+1)(k-s+2)+2$, then
$\rho\big(K_{n,n+s-k-1}-e\big)>\sqrt{\varepsilon_0(s)}$.
\\
{\em(ii)} If $n\geq 4k(k+1)$, then $\mu\big(K_{n,n+s-k-1}-e\big)>n+\frac{\varepsilon_0(s)}{n}.$
\end{Corollary}
\noindent{\bf  Proof:} In proofs below, denote $G=K_{n,n+s-k-1}-e$ and use $\rho$ and $\mu$ for $\rho(G)$ and $\mu(G)$,
respectively. Define
\begin{eqnarray*}
\Psi_1(\theta) & = & \theta^4 -\big(n(n+s-k-1)-1\big)\theta^2 +(n+s-k-2)(n - 1),\,\,\text{and}
\\
\Psi_2(\theta) & = & \theta^3-2(2n+s-k-2)\theta^2+\big((k-s)^2+(n - 1)(5n+5s-5k-6)\big)\theta\\\hspace{10pt}&&-(n - 1)(n+s-k-2)(2n+s-k-1).
\end{eqnarray*}
By setting $q=n+s-k-1$ and $\alpha\in \{0,1\}$ in  Lemma  \ref{33l},
$\rho$ and $\mu$ are equal to the maximum roots of $\Psi_1(\theta)$ and $\Psi_2(\theta)$, respectively.
As $n\geq (k+1)(k-s+2)+2>k-s+4$ and by $\Psi_1(\rho)=0$ it follows that
\begin{align*}
\rho^{2}&=\frac{1}{2}\left(n(n+s-k-1)-1+\sqrt{\big(n(n+s-k-1)-1\big)^2-4(n+s-k-2)(n - 1)}\,\,\right)\\&>n(n+s-k-1)-2\geq \varepsilon_0(s).
\end{align*}
This completes the proof of (i).

To prove (ii),  we first prove the following   claim.

\noindent
{\bf Claim 1. } $\Psi_2\left(2n+s-k-2+\frac{(k+1)(k+2)}{n}\right)<0.$
\par
By algebraic manipulations, we have
$$\Psi_2\left(2n+s-k-2+\frac{(k+1)(k+2)}{n}\right)=-\frac{1}{n^3}\Psi_3(n),$$ where $\Psi_3(n)= n^5-\Big(k(k+ 4)- s + 5\Big)n^4-\Big((k+2)(k+1)s-(k + 2)^3+s\Big)n^3-(k + 2)(k + 1)\big(2k^2 + 7k - s + 6\big)n^2+( k + 2)^2(k + 1)^2(k - s + 2)n -(k + 2)^3(k + 1)^3$.

Recall that  $n\geq 4k(k+1)$. Thus, $\Psi'''_3(n)=6\Big(10n^2-4\big(k(k+ 4)- s + 5\big)n -(k+2)(k+1)s+(k + 2)^3-s\Big)\geq \Psi'''_3(4k(k+1))=6\Big( s\big(15k^2 + 13k - 3\big)+k\big(144k^3 + 241k^2 + 22k - 68\big)+ 8\Big)>0$, and so  $\Psi''_3(n)\geq  \Psi''_3(4k(k+1))=2(k + 1)\Big(s\big(84k^3 + 60k^2 - 35k + 2\big)+k^2\big(544k^3 + 812k^2 - 154k - 347\big)+76k - 12\Big)>0$.

This leads to $\Psi'_3(n)\geq  \Psi'_3(4k(k+1))=(k+1)^2\Big(s\big(208k^4 + 112k^3 - 137k^2 + 12k - 4\big)+k\big(1024k^5 + 1328k^4 - 752k^3 - 791k^2 + 230k - 84\big)+8\Big)>0$, and so
$\Psi_3(n)\geq  \Psi_3(4k(k+1))=(k+1)^3\Big(4ks\big(48k^4 + 16k^3 - 45k^2 + 4k- 4\big)+k^2\big(768k^5 + 832k^4 - 928k^3 - 684k^2 + 215k - 150\big)+ 4(5k -2)\Big)>0$.
This completes the proof of Claim 1.

Direct computation yields that  $\Psi_2(0)=-(n - 1)(n+s-k-2)(2n+s-k-1)<0$ and  $\Psi_2(n+s-k-2)=n+s-k-2>0$. It is
observed that $\Psi_2(\theta)$ tends to infinity when $\theta$ tends to infinity.
Combining this with  $\Psi_2\left(2n+s-k-2+\frac{(k+1)(k+2)}{n}\right)<0$ by Claim 1, we   conclude that
$$\mu>2n+s-k-2+\frac{(k+1)(k+2)}{n}\geq n+\frac{\varepsilon_0(s)}{n},$$
and so  (ii) follows.  \q

\noindent
{\bf Proof of Theorem \ref{12t}.} Since $K_{n,n+s-k-1}-e\subset F^{0}_{n,k,s}$, by Corollary  \ref{31c},  $\rho\big(F^{0}_{n,k,s}\big)>\sqrt{\varepsilon_0(s)}$ and $\mu\big(F^{0}_{n,k,s}\big)>n+\frac{\varepsilon_0(s)}{n}$.
Thus Theorem \ref{12t} follows from Lemma \ref{32l} and Theorem \ref{11t}. \q

\section{The proofs of Theorem \ref{02c} and   Proposition \ref{31p}}

By Definition \ref{def} (vi) and (\ref{F=Z}), for an edge $w_0z_0 \in E(F_{n,k,s})$ with $d_{F_{n,k,s}}(w_0)=n+s-k-1$ and $d_{F_{n,k,s}}(z_0)=n$,
$F^{0}_{n,k,s}=F_{n,k,s}-w_0z_0$. {\bf
Throughout this section, we let $G= F^{0}_{n,k,s}$ and $G_0=F^{0}_{n,k,0}$.  Unless specially indicated,  let  $k$ and $s$ be two nonnegative integers such that $k\geq \max\{s,1\}$.}

\begin{Lemma} \label{rhoG} If   $n\geq (k+1)(k-s+2)+2$, then
$\rho_2(F^{0}_{n,k,s})<\sqrt{n(n+s-k-1)}$.
\end{Lemma}
\noindent {\bf Proof. }
By Corollary \ref{31c} and as $K_{n,n+s-k-1}-e\subset F^{0}_{n,k,s}=G$,
we have  $\rho(G)>\sqrt{\varepsilon_0(s)}$. By Lemma \ref{26l}, it follows that \begin{align*}\rho_2(G)&<\sqrt{|E(G)|-\varepsilon_0(s)}\\&=\sqrt{n^2-(k+1-s)(n-k)-1-n(n+s-k-2)-(k+1)(k+2-s)}
=\sqrt{n-( 2k- s + 3)}.\end{align*}
Since $n>2(k+1)$, $n(n+s-k-1)-\big(n-( 2k- s + 3)\big)=
n^2 -(k+2-s)n + 2k - s + 3\geq 2k(k+s)+4k + s+3>0$.
This completes the proof of the lemma.
\q.

\begin{Lemma} \label{62l} If $n\geq 2(k+1)-s$, then   $\rho(F^0_{n,k,s})$ is equal to the maximum root of $\Psi_4(x)$, where $\Psi_4(x)=x^2\Big(x^2-(n+s-k-2)\Big)\Big(x^2-k(k+1-s)\Big)
-\Big(x^2+(n+s-k-2)(x^2-1)\Big)\Big(kx^2+(n-k-1)(x^2-k(k+1-s))\Big).
$
\end{Lemma}
\noindent {\bf Proof.}  By Definition \ref{def} (vi) and (\ref{F=Z}), $G = [U,V]$ is a
bipartite graph and we may assume that $V$
contains $k+1-s$ vertices of degree $k$ and $U$ contains $k$ vertices of degree $n$ in $G$.
Define
$U_1 = \big\{u \in U: d_G(u) = n\big\}$ and $U_2=U\setminus U_1$,
$V_1 = \big\{v \in V: n-1 \le d_G(v) \le n\big\}$ and $V_2=V\setminus V_1$. By symmetry,
we may assume that $w_0\in U_2$ and $z_0\in V_1$.

Let  $f$ be the Perron vector of $G$, and let $\rho=\rho(G)$.
We shall adopt the following notation in the rest of the arguments:
\begin{align}\label{61e}\left\{
\begin{array}{llllll}
&x_1=f(w) \,\,\, \mbox{ if }\,\,\, w\in U_1,\\
& x_2=f(w)  \,\,\,\mbox{ if }\,\,\, w\in  V_2,
\\
&x_3=f(w)\,\,\,  \mbox{ if }\,\,\, w\in V_1\setminus \{z_0\},\\
&x_4=f(w)\,\,\, \mbox{ if }\,\,\, w\in U_2\setminus \{w_0\},\,\,
\\ &x_5=f(w_0),\,\,\text{and}\,\,x_6=f(z_0).\end{array}
\right.\end{align}  As $\big(A(G)\big)f=\rho f$, it follows that
\begin{align}\label{32e}
\left\{
\begin{array}{llllll}
\rho\,x_1=(k+1-s)x_2+(n+s-k-2)x_3+x_6, \\[1mm]
\rho \,x_2=kx_1, \\[1mm]
\rho\,x_3=kx_1+(n-k-1)x_4+x_5, \\[1mm]
\rho\,x_4=(n+s-k-2)x_3+x_6, \\[1mm]
\rho\,x_5=(n+s-k-2)x_3, \\[1mm]
\rho\,x_6=kx_1+(n-k-1)x_4.
\end{array}
\right.\end{align}
The first four equations of  (\ref{32e})  imply  that
\begin{align}
\label{33e}
x_4=\left(1-\frac{k (k+1-s)}{\rho^{2}}\right)x_1.
\end{align}
The equations on $x_3$, $x_5$ and $x_6$ of (\ref{32e}) lead to
\begin{align}\label{34e}
x_3= \frac{\rho^{2}}{\rho^{2}-(n+s-k-2)}x_6.
\end{align}
It follows from (\ref{34e}) and the first two equations of (\ref{32e}) that
\begin{align}\label{35e}
x_6= \frac{\big(\rho^{2}-(n+s-k-2)\big)\big(\rho^{2}-k(k+1-s)\big)}{\rho\big(\rho^{2}+(n+s-k-2)(\rho^2-1)\big)}x_1.
\end{align}
\par
With algebraic manipulations and utilizing (\ref{33e}), (\ref{35e}) and the sixth equation of (\ref{32e}),
$\rho$ is equal to  the maximum root of $\Psi_4(x)$, as desired.\q

\noindent {\bf Proof of Proposition \ref{31p}.} Let $\Psi_4(x)$ as defined in Lemma \ref{62l}.

To complete the proof , by Lemma \ref{rhoG},  it suffices to show that  \begin{align}
\label{36e}\Psi_4\left(\sqrt{n(n+s-k-1)}\right)>0.
\end{align}

Let $\Phi_2(n) = \Psi_4\Big(\sqrt{n(n+s-k-1)}\Big)$. Algebraic manipulation yields
$\Phi_2(n) = 2n^4-\Big((k^2 + 4)(k+1 - s)+2\Big)n^3+\Big(2k^3(k - 2s + 2)
+2k^2(s - 1)^2-(2s-5)(k-s)-2(s -3)\Big)n^2-
(k - s + 1)(k - s + 2)\Big(k^2(k- s)- 2k + 1\Big)n -k(k + 1)(k - s + 1)(k - s + 2)$.

\noindent
{\bf Case 1.} $s\geq 1$.

Let $\Phi_3(n)=\Big(2k^3(k - 2s + 2)+2k^2(s - 1)^2-(2s-5)(k-s)-2(s -3)\Big)n^2-
(k - s + 1)(k - s + 2)\Big(k^2(k- s)- 2k + 1\Big)n -k(k + 1)(k - s + 1)(k - s + 2)$ and
$\Phi_4(n)= 2n^4-\big((k^2 + 4)(k+1 - s)+2\big)n^3$.

When $2n\geq(k^2+4)(k+1)$, we have $\Phi'_3(n)\geq \Phi'_3\Big(\frac{1}{2}(k^2+4)(k+1)\Big)=
2k^5(k - s)^2+2k^4(k-s)(3k-s)+k^3(k - s)(13k - 6s+28)+k^3(s^2-7s+13)+7k^2(s-1)s+
22k^2(k-s)+k^2(s^3-4s+4) + 5k(2s^2 - 8s + 7k) + 45k-25s + 7s^2+ 22>0$. This implies that
$\Phi_3(n)\geq \Phi_3\Big(\frac{1}{2}(k^2+4)(k+1)\Big)=\frac{1}{4}(k+1)\Big(2k^7(k-s)^2+2k^6(k-s)(3k-s)+4k^5(k-s)(5k-3s)+k^4(k-s)(49k-12s)
+k^3(k-s)(83k - 16s)+12k^3\big(s(s-1)+12(k-s)\big)+2k^2s\big(k^2s^2 - 2k + 4s^2\big)+2k^2(11s^2 - 77s + 91k)+4(46k^2 - 41ks + 6s^2)+44(ks^2 - 2s + 4k)+80\Big)>0.$

When $2n\geq(k^2+4)(k+1)$, we have  $2n-\big((k^2 + 4)(k+1 - s)+2\big)\geq s(k^2 + 4) - 2\geq k^2+2>0$, and so
$\Phi_4(n)>0$.
As $\Phi_2(n)=\Phi_3(n)+\Phi_4(n)>0$, it follows that  (\ref{36e}) must hold.

\noindent
{\bf Case 2.} $s=0$.

Define $\Phi_5(n)=2n^3 -\Big(k^2(k+1)+2k+4\Big)n^2 + (k + 2)\big(k^3 - k + 1\big)n + k(k + 2)(k + 1)$.
As $s=0$, $\Phi_2(n)=(n-k-1)\Phi_5(n)$.
Since $2n\geq(k^2+4)(k+1)>k^2(k+1)+2k+4$, we have $\Phi_5(n)>(k + 2)\big(k^3 - k + 1\big)n + k(k + 2)(k + 1)> k(k + 2)(k + 1)>0$.
Thus, $\Phi_2(n)>0$ and so (\ref{36e}) holds.
\q

\begin{Lemma} \label{63l} If $n\geq 3k(k+1)$ and $k\geq 2$, then   $\rho(F^{0}_{n,k,0})>\sqrt{n(n-k-2)+(k+2)^2}$.
\end{Lemma}
\noindent {\bf Proof.} Throughout this proof, we  simplify  rewrite   $\rho(G_0)$  as $\rho$.

 By Lemma \ref{62l}, $\rho$ is equal to the maximum root of  $\Psi_5(x)$, where $\Psi_5(x)=x^2\Big(x^2-(n-k-2)\Big)\Big(x^2-k(k+1)\Big)
-\Big(x^2+(n-k-2)(x^2-1)\Big)\Big(kx^2+(n-k-1)(x^2-k(k+1))\Big).$ To show that $\rho>\sqrt{n(n-k-2)+(k+2)^2}$, it suffices to prove $\Psi_5(\sqrt{n(n-k-2)+(k+2)^2})<0.$

Denote by  $\Phi_6(n)=\Phi(\sqrt{n(n-k-2)+(k+2)^2})=- n^5+(k^2 + 6k + 10)n^4-\big(3k(k^2 + 5k + 13)+37\big)n^3+(k + 2)(4k^3 + 17k^2 + 45k + 44)n^2 -(k^3(3k^2 + 22k + 78)+174k^2 + 219k + 112)n + (k + 2)(k^5 + 7k^4 + 26k^3 + 64k^2 + 85k + 44)$.

When  $n\geq 3k(k+1)$ and $k\geq 2$, we have $\Phi{''''}_6(n)=24\big(-5n+k^2+6k+10\big)\leq -24(14k^2 + 9k - 10)<0$ and thus  $\Phi{'''}_6(n)\leq \Phi{'''}_6(3k(k+1))= -6\big(k^2(78k^2-87)+9k(11k^2-9)+ 37\big)<0$. Once again, since $n\geq 3k(k+1)$ and $k\geq 2$, we have  $\Phi{''}_6(n)\leq \Phi{''}_6(3k(k+1))=-2\big(k^4(216k^2-274) +k^3(405k^2 - 673)+ 65k^2 + 199k - 88\big)<0$ and so
$\Phi'_6(n)\leq \Phi'_6(3k(k+1))=-\big(k^3(297k^5 + 729k^4- 375k^3- 1899k^2- 575k+303)+ 159k^2(k-1)+309k(k^2-1)+ 112\big)<0.$

 This implies that $\Phi_6(n)\leq \Phi_6(3k(k+1))=
 -\big(162k^{10} + 486k^9 - 198k^8 - 1827k^7 - 1339k^6 + 516k^5 + 41k^4 - 728k^3 - 12k^2 + 122k - 88\big)<-k^3\big(162k^7+ 486k^6- 198k^5- 1827k^4- 1339k^3+ 516k^2- 728\big)<-k^6\big(162k^4+ 486k^3- 198k^2- 1827k- 1339\big)$. Denote by $\Phi_7(k)=162k^4+ 486k^3- 198k^2- 1827k- 1339$. Since $\Phi_7(2)=695>0$ and  $\Phi_7(k)=175k^3-1339+k^2(311k-198)+k(162k^3- 1827)>0$ for $k\geq 3$, we have $\Phi_6(n)<0$, as desired. \q

 \begin{Lemma} \label{64l} If $n\geq 2(k+2)^2$ and $k\geq 1$, then   $\mu(F^{0}_{n,k,0})>2n-k-1.5$.
\end{Lemma}
\noindent {\bf Proof.} Throughout this proof, we  simplify  rewrite   $\mu(G_0)$ as $\mu$. By replacing $s=0$  in the proof of Lemma \ref{62l}, we define  $U_1$, $U_2$,
$V_1$ and $V_2$, where  $w_0\in U_2$ and $z_0\in V_1$. Moreover, we let $f$ be the Perron vector of $G_0$ and we  also adopt the same  notation from (\ref{61e}). As $\big(Q(G_0)\big)f=\mu f$, it follows that
\begin{align}\label{62e}
\left\{
\begin{array}{llllll}
(\mu-n)\,x_1=(k+1)x_2+(n-k-2)x_3+x_6, \\[1mm]
(\mu-k)\,x_2=kx_1, \\[1mm]
(\mu-n)\,x_3=kx_1+(n-k-1)x_4+x_5, \\[1mm]
(\mu-n+k+1)\,x_4=(n-k-2)x_3+x_6, \\[1mm]
(\mu-n+k+2)\,x_5=(n-k-2)x_3, \\[1mm]
(\mu-n+1)x_6=kx_1+(n-k-1)x_4.
\end{array}
\right.\end{align}

The first four equations of  (\ref{62e})  imply  that
\begin{align}
\label{63e}
x_4=\frac{(\mu-n)(\mu-k)-k(k+1)}{(\mu-k)(\mu-n+k+1)}x_1.
\end{align}
The equations on $x_3$, $x_5$ and $x_6$ of (\ref{62e}) lead to
\begin{align}\label{64e}
x_3= \frac{(\mu-n+1)(\mu-n+k+2)}{(\mu-n)(\mu-n+k+2)-(n-k-2)}x_6.
\end{align}
It follows from (\ref{64e}) and the first two equations of (\ref{62e}) that
\begin{align}\label{65e}
 x_6= \frac{(\mu^2-( k + n)\mu+k(n-k-1))(\mu^2 -(2n-k-2)\mu+(n - 1)(n-k-2))}
{(\mu-k)((n - k - 1)\mu^2 -(2n^2 -(3k+5)n+(k+2)^2)x +(n - 1)(n-k -1)(n-k-2))}x_1.
\end{align}
\par
With algebraic manipulations and utilizing (\ref{63e}), (\ref{65e}) and the sixth equation of (\ref{62e}),
$\mu$ is equal to  the maximum root of $\Psi_6(x)$, where
$\Psi_6(x)=x^4 -\big(4n-4-k\big)x^3+\big(5n^2-( k + 11)n - 2k^2 + 6\big)x^2-\big(2n^3 + (2k - 7)n^2 -( 6k^2 + 7k - 7)n +2k^3 + 8k^2 + 6k - 2\big)x + 2 k(n - 1)(n-k-1)(n-k -  2).$

Since $\Psi_6(2n-k-1.5)=\frac{-1}{16}\Phi_8(n),$ where $\Phi_8(n)=16n^3+4(4k^2 - 4k - 19)n^2 -4(4k^3 - 25k - 23)n - 8k^3 - 36k^2 - 66k - 33=4n(4n^2-19n-4k^3+25k+23)+16k(k-1)n^2- 8k^3 - 36k^2 - 66k - 33$. Note that $n\geq 2(k+2)^2$. Thus,  $\Phi''_8(n)=8( 12n +4k^2 - 4k- 19)\geq \Phi''_8(2(k+2)^2)=8(28k^2 + 92k + 77)>0$ and so $\Phi'_8(n)\geq \Phi'_8(2(k+2)^2)= 4(64k^4 + 428k^3 + 1076k^2 + 1193k + 487)>0$. This implies that $\Phi_8(n)\geq \Phi_8(2(k+2)^2)=192k^6 + 1952k^5 + 8272k^4 + 18624k^3 + 23348k^2 + 15294k + 4031$. Now, we can conclude that $\Psi_6(2n-k-1.5)<0$,  completing   the proof of this result.\q

\begin{Lemma}\label{81l}{\em (Li and Ning \cite{Ning2017})} Let $G_1$ be a balanced bipartite graph on $2n$ vertices. If $\delta(G_1)\geq k\geq 1$, $n\geq 2k+3$ and $|E(G_1)|>n(n-k-2)+(k+2)^2$, then $G_1$ is traceable  unless $G_1\subseteq F_{n,k,0}$ or $k=1$ and $G_1\subseteq K_{n-1,n-1}\cup K_2.$\end{Lemma}
\noindent
{\bf Proof of Theorem \ref{02c}.}  From    (\ref{22e}) and Lemmas \ref{63l}--\ref{64l},  we have
\begin{align*}
&n+\frac{n(n-k-2)+(k+2)^2}{n}\leq 2n-k-1.5<\mu(G_0) \leq  \frac{|E(G_0)|}{n}+n,\\&\text{and}\,\,\sqrt{n(n-k-2)+(k+2)^2}< \rho(G_0) \leq  \sqrt{|E(G_0)|},
\end{align*}
which implies that $|E(G_0)|>n(n-k-2)+(k+2)^2$.   From Lemma \ref{81l} and since $\delta(G_0)\geq 2$, we have   $G_0\subseteq F_{n,k,0}$. Now, the result    follows from Lemma \ref{32l}.\q
\section{The proof  of Theorem  \ref{13t}}

Throughout this section, we assume that
$G=[U,V]$ is an   almost  balanced bipartite graph with  $|U|=|V|+1=n$.
Let $v_0$ be a vertex not in $V(G)$ and define a balanced bipartite graph $G^{v_0}$
from $G$ by adding $v_0$ and $n$ edges joining $v_0$ to all vertices of $U$.

\begin{Lemma}\label{11o}
If  $G^{v_0}$ is $(q,q)$-Hamiltonian, then  $G$  is  $(q,q)$-traceable.
\end{Lemma}
\noindent{\bf  Proof:} As the case when $q=0$ follows from definition immediately,
we assume that $q\geq 1$.    Let $S$ be an arbitrary set of $2q$ vertices of $G$ such that $|S\cap U|=q=|S\cap V|$.
Choose a vertex $v \in S \cap V$. Let $S_1=\big(S\setminus \{v\}\big)\cup \{v_0\}$ and $V_1=V\cup \{v_0\}$.
Then $|S_1\cap U|=q=|S_1\cap V_1|$. Since $G^{v_0}$ is $(q,q)$-Hamiltonian,
$G^{v_0}\big[V\left(G^{v_0}\right)\setminus S_1\big]$ contains a Hamiltonian cycle,
and hence $G\big[V(G)\setminus S\big]$ is traceable, as
$G^{v_0}\big[\left(V\left(G^{v_0}\right)\setminus S_1\right)\setminus\{v\}\big]= G\big[V(G)\setminus S\big]$.
By the arbitrariness  of $S$, $G$  is  $(q,q)$-traceable.
\q

\noindent{\bf  Proof of Theorem \ref{13t} (i):}
We first show that, under the assumption of Theorem \ref{13t}, we have
\begin{equation}
\label{eG}
|E(G)| >n(n+q-k-2)+(k+1)(k+1-q).
\end{equation}
Assume that $|E(G)| \le n(n+q-k-2)+(k+1)(k+1-q)$.
By  (\ref{omega}) and Corollary \ref{C4},
we have  $\Omega(\alpha) \ge \alpha \left(\frac{|E(G)|}{n}+n\right)+(1-\alpha)\sqrt{|E(G)|} \ge \Theta(G,\alpha)$,
contrary to the assumption that $\Omega(\alpha)<\Theta(G,\alpha)$. Hence (\ref{eG}) follows.

From (\ref{eG}), it follows that
$|E(G^{v_0})| >n(n+q-k-1)+(k+1)(k+1-q)=\varepsilon_0(q+ 1)$. By Theorem \ref{21t} and Proposition \ref{11p},
either  $G^{v_0}$ is   $(q,q)$-Hamiltonian  or $G^{v_0}\subseteq F_{n,k,q+1}$.
It follows by Lemma  \ref{11o} that either  $G$ is $(q,q)$-traceable  or $G\subseteq Z_{n+q-k-1,k-q}$.
\q

\noindent{\bf  Proof of Theorem \ref{13t} (ii):}
By Corollary \ref{31c}, we have $\rho\big(K_{n,n+q-k-1}-e\big)>\sqrt{\varepsilon_0(q)}>\Omega(0)$ and $\mu\big(K_{n,n+q-k-1}-e\big)>n+\frac{\varepsilon_0(q)}{n}>\Omega(1)$. Note that  $K_{n,n+q-k-1}-e\subset Z^{0}_{n+q-k-1,k-q}$.
 Thus, Theorem \ref{13t} (i) implies that $G$ is $(q,q)$-traceable  unless  $G\subseteq Z_{n+q-k-1,k-q}$.
 Now, the result follows from Lemma \ref{32l}.
\q
\par\medskip
\noindent {\bf Acknowledgment.}   The authors would like  to thank  three anonymous  referees
for their valuable comments which lead to a great  improvement of the
original manuscript.

 \end{document}